\documentclass[12pt]{amsart} 
\usepackage{amssymb}

\def\page{\marginpar}
\def\marginpar#1{\ignorespaces}

\makeatletter

\usepackage{euscript}

\let\cal\EuScript

\def\cases#1{\left\{\,\vcenter{\normalbaselines\m@th
    \ialign{$##\hfil$&\quad##\hfil\crcr#1\crcr}}\right.}
\def\matrix#1{\null\,\vcenter{\normalbaselines\m@th
    \ialign{\hfil$##$\hfil&&\quad\hfil$##$\hfil\crcr
      \mathstrut\crcr\noalign{\kern-\baselineskip}
      #1\crcr\mathstrut\crcr\noalign{\kern-\baselineskip}}}\,}

\newtheorem{theorem}{Theorem}[section]

\newtheorem{corollary}[theorem]{Corollary}
\newtheorem{proposition}[theorem]{Proposition}

\theoremstyle{definition}

\def\eqlabel#1{\label{eq#1}}
\def\eqref#1{(\ref{eq#1})}
\numberwithin{equation}{section}

\let\<\langle
\let\>\rangle

\newcommand\Z{\mathbb Z}

\newcommand\R{\mathbb R}

\newcommand\PP{\mathbb P}
\newcommand\N{\mathbb N}
\newcommand\E{\mathbb E}

\def\exp{\operatorname{exp}}

\let\over\@@over
\let\atop\@@atop
\let\above\@@above
\let\overwithdelims\@@overwithdelims
\let\atopwithdelims\@@atopwithdelims
\let\abovewithdelims\@@abovewithdelims

\def\eqalign#1{\null\,\vcenter{\openup\jot
  \ialign{\strut\hfil$\displaystyle{##}$&$\displaystyle{{}##}$\hfil
      \crcr#1\crcr}}\,}
\def\eqalignbot#1{\null\,\vbox{\openup\jot
  \ialign{\strut\hfil$\displaystyle{##}$&$\displaystyle{{}##}$\hfil
      \crcr#1\crcr}}\,}

\makeatother

\begin{document}

\title[Smoluchowski's coagulation equation]{Smoluchowski's coagulation equation:~uniqueness, non-uniqueness
and a hydrodynamic limit for the stochastic coalescent}
\author{James R. Norris}
\address{Statistical Laboratory\\
16 Mill Lane\\
Cambridge\\
CB2 15B\\
United Kingdom}
\email{j.r.norris@statslab.cam.ac.uk}
%\urladdr{www.statslab.cam.ac.uk/{\tt\char`\~}james}
%\curraddr{MSRI\\ 1000 Centennial Drive\\ Berkeley, CA  94720-5070\\ USA}

\begin{abstract}
Sufficient conditions are given for existence and uniqueness in
Smoluchowski's coagulation equation, for a wide class of coagulation
kernels and initial mass distributions. An example of non-uniqueness is
constructed. The stochastic coalescent is shown to converge weakly to
the solution of Smoluchowski's equation.
\end{abstract}

\thanks{This research was supported financially by the European Union under
contract FMRX CT96 0075A and
by the Mathematical Sciences Research Institute.
Research at MSRI is supported in part by NSF grant DMS-9701755.}
\maketitle

\section{Introduction}
\page{2}
Coagulation of particles, in pairs and over time, within a large
system of particles, is a phenomenon widely observed and widely
postulated in scientific models. Examples arise in the study of
aerosols, of phase separation in liquid mixtures, in polymerization and
astronomy. Typically, it is argued that the rate at which pairs of
particles coagulate depends, for physical reasons, in a given way, on
some positive parameter associated to each particle, such as size or
mass. In the models we shall consider, it is further argued that the
effects of spatial fluctuations in the size or mass density are
negligible---for example, by supposing that the particles perform
independent random motions on a time scale faster than the process of
coagulation. 
\page{3}

The first mathematical model of this sort of process was proposed by
Smoluchowski \cite{Sm} in 1916, see also  Chandrasekhar \cite{Ch}.
Smoluchowski argued that particles of radius $r$ would perform
independent Brownian motions of variance proportional to $ 1/r $, so
pairs of particles of radii $ r_1 $ and $ r_2 $ would meet at a rate
proportional to
$$
(r_1 + r_2) (1/r_1 + 1/r_2).
$$
\marginpar{Added period}
Expressed in terms of masses, this leads to the coagulation kernel
$$
K(x, y) = (x^{1/3} + y^{1/3})(x^{-1/3} + y^{-1/3})
$$
for particles of masses $x$ and $y$. Then, making some implicit
assumptions about ergodic averages, Smoluchowski wrote down the
following infinite system of differential equations for the evolution
of densities $ \mu(x) $ of particles of mass $ x = 1,2,3, \ldots $
$$
\frac{d}{dt} \mu_t (x) = \frac{1}{2} \sum^{x - 1}_{y = 1} K(y, x - y) \mu_t (y)
\mu_t (x - y) - \mu_t (x) \sum^\infty_{y = 1} K(x,y) \mu_t (y).
$$
\page{4}
Here, the first sum on the right corresponds to coagulation of smaller
particles to produce one of mass $x$, whereas the second sum
corresponds to removal of particles of mass $x$ as they in turn
coagulate to produce larger particles.

More generally, in other models, such systems of equations are
considered for many different coagulation kernels $K$, see for
example \cite{A}. Also, analogous integro-differential equations are
considered which allow for a continuum of masses $x$. It is known by
now that, for a suitable initial mass distribution $ \mu_0 $,
Smoluchowski's original equations have a unique solution. Much progress has been
made in determining when existence and uniqueness holds for more
general coagulation kernels, see \cite{McL1}, \cite{McL2}, \cite{W},
\cite{BC}, \cite{H} for discrete mass distributions. Nevertheless
many fundamental questions remain open, even for certain coagulation kernels 
studied
extensively in applied sciences.
\page{5}

In this paper we give some new positive results on the existence and
uniqueness problem for Smoluchowski-type equations. In particular:
\begin{itemize}
\item we prove existence of solutions for continuous coagulation
kernels $K$ such that
$$
K(x,y)/xy \rightarrow 0  \quad \text{ as } (x,y) \rightarrow \infty 
$$
extending a result of Jeon \cite{Je} for the discrete case;
\item we prove local existence and uniqueness of solutions when $ K(x,y)
\leq \varphi(x)\varphi(y) $ for some continuous sublinear function $ \varphi : E
\rightarrow  (0, \infty) $, provided that the initial mass distribution $
\mu_0 $ satisfies
$$
\int_{(0, \infty)} \varphi(x)^2 \mu_0 \, (dx)<\infty;
$$
\item this allows us to treat the case where $ K(x,y) $ blows up as $
x \rightarrow 0 $ or $ y \rightarrow 0 $, also to prove uniqueness in some cases
when the mass distribution has no second, or even first, moment;
\item we can do without any local regularity conditions on $K$;
\item we do not have to assume that the initial mass distribution is
discrete, nor that it has a density with respect to Lebesgue measure.
\end{itemize}
\page{6}

We also construct in \S3 an example of a coagulation kernel $K$ and an
initial mass distribution $ \mu_0 $, such that Smoluchowski's equation
has at least two distinct solutions, both of which are conservative,
in the sense that
$$
\int_{(0,\infty)} x \mu_t \, (dx) = \int_{(0, \infty)} x \mu_0 \, (dx)<\infty
$$
for all $t$.

Then in \S4 we consider a stochastic system of coagulating particles,
\marginpar{Check running\\title above}
where particles of masses $x$ and $y$ coagulate at a rate proportional
to $ K(x,y) $. We show that, when we can establish uniqueness in
Smoluchowski's equation, the particle system, suitably normalized,
converges weakly to the solution of the deterministic equation. Thus
we obtain a statistical derivation of Smoluchowski's equation.
\page{7}
This goes some way towards resolving Problem 10 in Aldous' survey article \cite{A}.

\page{8}
\section{Existence and uniqueness in Smoluchowski's coagulation equation}
\label{2}
Let $ E = (0, \infty) $ and let $ K: E \times E \rightarrow [0,
\infty) $ be a symmetric measurable function, the {\em  coagulation
kernel}. Denote by $\cal M=\cal M_E$
the space of signed Radon measures on $E$, that is to say, those
signed measures having finite total variation on each compact subset
of $E$. Denote by $ \cal M^+ $ the set of (non-negative) measures in $\cal M$.
If $ \mu \in \cal M^+ $ satisfies, for all compact sets $ B \subseteq E$
$$
\int_{B \times E} K (x,y) \mu \, (dx) \mu (dy) < \infty,
$$
then we define $ L(\mu) \in \cal M $ by
$$
\<f, L(\mu) \> = \frac{1}{2} \int_{E \times E} \{ f(x + y) - f (x) - f(y) \}
K(x,y) \mu(dx) \mu (dy)
$$
for all bounded measurable functions $f$ of compact support.
\page{9}

We consider  the following weak form of Smoluchowski's coagulation
equation
\begin{equation}
\eqlabel{C}
\mu_t = \mu_0 + \int^t_0 L(\mu_s) \, ds.
\end{equation}
We admit as a {\em  local solution} any map
$$
t \mapsto \mu_t: [0,T) \mapsto \cal M^+
$$
where $ T \in (0, \infty] $, such that:
\begin{itemize}
\item[(i)] we have
\marginpar{moved (i) here}
$$
\int_E x 1_{x \leq 1} \,\mu_0 (dx) < \infty;
$$
\item[(ii)] for all compact sets $ B \subseteq E $, the following map is
measurable
$$
t \mapsto \mu_t(B): [0,T) \rightarrow [0, \infty);
$$
\item[(iii)] we have, for all $ t < T $ and all compact sets $ B
\subseteq E $
$$
\int^t_0 \int_{B \times E} K (x,y) \mu_s (dx) \mu_s (dy) ds < \infty;
$$
\item[(iv)] for all bounded measurable functions $f$ of compact
support and also for $ f(x) = x 1_{ x \leq 1} $, for all $ t < T $
\begin{equation}
\eqlabel{fL}
\<f, \mu_t \> = \<f, \mu_0 \> + \int^t_0 \<f, L(\mu_s)\> ds.
\end{equation}
\end{itemize}
In the case $ T = \infty $, we have a {\em  solution}.
\page{10}

The condition that, for $ f(x) = x1_{x \leq 1} $, we have $ \<f,
\mu_0 \>  < \infty $ and that \eqref{fL} holds 
\marginpar{removed comma}
is a boundary
condition, expressing that no mass enters at $0$. We obtain an
equivalent condition on replacing $f$ by any non-vanishing sublinear function 
$E\rightarrow [0,\infty)$ of
bounded support, which is linear near $0 $. 
A function $f:E\rightarrow [0,\infty)$ is {\it sublinear} if 
$$
f(\lambda x)\leq\lambda f(x)\quad\text{ for all }\quad x\in E, \lambda\geq 1.
$$
Note that such a function $f$ is always {\it subadditive}:
$$
f(x+y)\leq f(x)+f(y)\quad\text{ for all }\quad x,y\in E.  
$$
Hence $\< f,L(\mu)\>\leq 0$ for all $\mu\in\cal M^+$.
Note also that, if $
\varphi: E \rightarrow [0, \infty) $ is {\it any} sublinear function and if
$$
\varphi_n (x) = \cases{
n x \varphi(n^{-1}), & $ 0 < x \leq n^{-1}$,\cr
\varphi(x), & $ n^{-1} < x \leq n $,\cr
0, & $ x >  n $,
}
$$
then $ \varphi_n (x) \uparrow \varphi(x) $ for all $x$, and $
\varphi_n $ is sublinear of bounded support, linear near $ 0 $. So,
for $ t < T $ 
$$ 
\< \varphi_n, \mu_t\> - \<\varphi_n, \mu_0\> = \int^t_0 \<\varphi_n,
L(\mu_s)\> \leq 0.
$$
\page{11}
Hence, using monotone convergence on the left and Fatou's lemma on the
right
\begin{equation}
\eqlabel{F}
\<\varphi, \mu_0 \> \geq \<\varphi, \mu_t \> - \int^t_0 \<\varphi,
L(\mu_s)\> \, d s.
\end{equation}
In particular, $ \<\varphi, \mu_t\> $ is non-increasing in $t$. In
particular, the {\em  total mass density}
$$
\int_E x \mu_t \, (dx)
$$
is non-increasing in $t$; if it is finite and constant, we say that $
(\mu_t)_{t < T} $ is {\em  conservative}.
\page{12}

Throughout this section we make the basic assumption that
\begin{equation}
\eqlabel{varphiK}
K(x,y) \leq \varphi(x) \varphi(y) \quad \text{ for all } x, y \in E
\end{equation}
where $ \varphi: E \rightarrow (0, \infty) $ is a continuous sublinear
function. We also assume that the initial measure $ \mu_0 $ satisfies
\begin{equation}
\eqlabel{varphimu}
\<\varphi, \mu_0\> < \infty.
\end{equation}
We call any local solution $ (\mu_t)_{t < T} $ such that 
$$
\int^t_0  \<\varphi^2, \mu_s \> \, ds < \infty \quad \text{ for all }
t < T
$$
a {\em  strong} local solution.
\page{12a}

Here is a summary of what is known so far about existence, uniqueness
and conservation of mass in Smoluchowski's equation. The picture is
more complete for discrete mass distributions---that is when $ \mu_0 $
is supported on $ \N $. Then, provided $ \mu_0 $ has a finite second
moment, Ball and Carr \cite{BC} proved existence and mass conservation
when $ K (x,y) \leq x + y $, Heilman \cite{H} added uniqueness under
the same hypotheses. Jeon \cite{Je} has recently proved global
existence when $ K (x,y)/xy \rightarrow 0 $ as $ (x,y) \rightarrow
\infty $. McLeod \cite{McL1} long ago proved local existence when $ K
(x,y) \leq xy $. For general mass distributions $ \mu_0 $, less is
known. Dubovskii and Stewart \cite{DS} have shown existence and
uniqueness provided $ \mu_0 $ has an exponential moment and a
continuous density with respect to Lebesgue measure, and provided
$K$ is continuous with $ K (x,y) \leq 1 + x + y $. Recently, Clark and
Katsouros \cite{CK} proved existence and uniqueness for a particular choice of
kernel which blows up when $x$ or $y$ is small.
\page{13}
\begin{theorem}
\label{2.1}
Assume conditions  {\rm\eqref{varphiK}} and  {\rm\eqref{varphimu}}. If $
(\mu_t)_{t < T} $ and $ (\nu_t)_{t < T} $ are local solutions, starting
from $ \mu_0 $, and if $ (\mu_t)_{t < T} $ is strong, then $ \mu_t =
\nu_t $ for all $ t < T $. If $ \varphi(x) \geq  \varepsilon x $ for all
$x$, for some $ \varepsilon >  0 $, then any strong solution is
conservative. Moreover, if $ \<\varphi^2, \mu_0\> < \infty $, then
\begin{itemize}
\item[(i)] there exists a unique maximal strong solution $
(\mu_t)_{t < \zeta(\mu_0)},$ with $ \zeta(\mu_0) \geq \<\varphi^2,
\mu_0\>^{-1} , $
\item[(ii)] if $ \varphi^2 $ is sublinear or if $ K (x,y) \leq
\varphi(x) + \varphi(y) $ for all $ x,y \in E $, then $ \zeta (\mu_0)
= \infty $.
\end{itemize}
\end{theorem}
The proof will occupy the remainder of this section.
The method is to find an approximation to Smoluchowski's equation
by a system depending on $K$ and $\varphi$ only through their values on
a given compact set.
The idea of the approximation may be readily understood in terms
of the stochastic coalescent, for which a directly analogous 
approximation is discussed in \S 4.
Moreover, the finite particle interpretation explains certain crucial
non-negativity statements, which are given less transparent analytical proofs
below.

We remark that this theorem provides examples where uniqueness holds,
even when the solution fails to be conservative, in the trivial sense
that the total mass density is infinite. We have not yet found an example
of a strong solution which has finite initial mass density and which
fails to conserve mass. The example of \S 3 shows, on the other hand,
that uniqueness can fail while solutions remain conservative.

%\begin{proof}
Let $ B \subseteq E $ be compact. We will eventually pass to the limit
$ B \uparrow E $.
Denote by $\cal M_B $ the space of finite signed measures supported on
$B$. Note that $ \varphi $ is bounded on $B$. We define $ L^B: \cal M_B
\times \R \rightarrow \cal M_B \times \R $ by the requirement
\begin{multline*}
\<(f,a), L^B(\mu, \lambda)\> =\\
\frac{1}{2} \int_{E \times E}\{ f(x + y) 1_{x + y \in B} + a \varphi(x + y)
1_{x + y \not\in B} -f(x) - f(y) \}K (x,y) \mu (dx) \mu(dy) \\
+ \lambda \int_E \{ a \varphi(x) - f(x) \} \varphi(x) \mu(dx)
\end{multline*}
for all bounded measurable functions $f$ on $E$ and all $ a \in \R $.
Here we used the notation $\<(f,a),(\mu,\lambda)\>=\<f,\mu\>+a\lambda$.
\page{15}

Consider the equation
\begin{equation}
\eqlabel{B}
(\mu_t, \lambda_t) = (\mu_0, \lambda_0) + \int^t_0 L^B (\mu_s,
\lambda_s) \, ds.
\end{equation}
We admit as a {\em  local solution} any continuous map
$$
t \mapsto (\mu_t, \lambda_t ): [0,T] \rightarrow\cal M_B \times \R
$$
where $ T \in(0, \infty) $, which satisfies the equation for all $ t
\in [0,T] $. When $ [0,T] $ is replaced by $ [0, \infty) $ we get the
notion of {\em solution}.

\begin{proposition} 
\label{2.2}
Suppose $ \mu_0 \in\cal M_B $ with $ \mu_0 \geq 0 $ and that $ \lambda_0 \in
[0, \infty) $. The equation {\rm \eqref{B}} has a unique solution $ (\mu_t,
\lambda_t)_{t \geq 0} $ starting from $ (\mu_0, \lambda_0) $. Moreover
$ \mu_t \geq 0 $ and $ \lambda_t \geq 0 $ for all $t$.
\end{proposition}
\begin{proof}
Our basic assumption \eqref{varphiK} remains valid when $ \varphi $ is
replaced by $ \varphi + 1 $, so we may assume without loss that $
\varphi \geq 1 $. By a scaling argument we may assume, also without loss,
that
$$
\< \varphi, \mu_0 \> + \lambda_0 \leq 1
$$
which implies that
$$
\|\mu_0 \| + |\lambda_0| \leq 1.
$$
We shall show, by a standard iterative scheme, that there is a
constant $ T >  0 $, depending only on $ \varphi $ and $B$, and a
unique local solution $ (\mu_t, \lambda_t)_{t \leq T} $ starting from
$ (\mu_0, \lambda_0) $. Then we shall show, moreover, that $ \mu_t
\geq  0 $ for all $ t \in [0,T] $. 
\page{17}

First of all, let us see that this is enough to prove the proposition.
If we put $ f = 0 $ and $ a = 1 $ in \eqref{B}, we obtain
$$
\frac{d}{dt} \lambda_t = \frac{1}{2} \int_{E \times E} \varphi(x + y) 1_{x + y
\not\in B} K (x,y) \mu_t (dx) \mu_t (dy) + \lambda_t \int_E \varphi(x)^2 \mu_t (dx).
$$
So, since $ \mu_t \geq 0 $, we deduce $ \lambda_t \geq 0 $ for all
$t$. Next, we put $ f = \varphi $ and $ a = 1 $ to see that 
$$
\frac{d}{dt}( \< \varphi, \mu_t \> + \lambda_t) 
= \frac{1}{2} \int_{E \times E} \{ \varphi(x + y) - \varphi(x) - \varphi(y)
\} K (x,y) \mu_t (dx) \mu_t (dy) \leq 0.
$$
Hence
$$
\|\mu_T\| + |\lambda_T| \leq \< \varphi, \mu_T \> + \lambda_T 
 \leq \< \varphi, \mu_0 \> + \lambda_0 \leq 1.
$$
We can now start again from $ (\mu_T , \lambda_T) $ at time $T$ to
extend the solution to $ [0, 2T] $, and so on, to prove the
proposition.
\page{18}

We use the following norm on $ \cal M_B \times \R $:
$$
\|(\mu, \lambda)\| = \| \mu \| + |\lambda|.
$$
We note the following estimates:~there is a constant $ C < \infty $,
depending only on $ \varphi $ and $B$ such that, for all $ \mu, \mu '
\in\cal M_B $  and all $ \lambda, \lambda ' \in \R $,
\marginpar{Silvio: is this OK?}
\begin{equation}
\eqlabel{C1}
\|L^B(\mu, \lambda)\|  \leq C \|(\mu, \lambda)\|^2,
\end{equation}
\begin{equation}
\eqlabel{C2}
\|L^B(\mu, \lambda) - L^B(\mu ', \lambda ')\| 
\leq C \|(\mu, \lambda) - (\mu ', \lambda ') \|\, (\|(\mu, \lambda)\| +
\|(\mu ', \lambda ')\|).
\end{equation}
\page{19}
We turn to the iterative scheme. Set $ (\mu^0_t, \lambda^0_t) =
(\mu_0, \lambda_0 ) $ for all $t$ and define inductively a sequence of
continuous maps
$$
t \mapsto (\mu^n_t, \lambda^n_t): [0, \infty) \rightarrow\cal M_B \times
\R
$$
by
$$
(\mu^{n + 1}_t, \lambda^{n + 1}_t) = (\mu_0, \lambda_0) + \int^t_0 L^B
(\mu^n_s, \lambda^n_s) \, ds.
$$
Set
$$
f_n(t) = \|(\mu^n_t, \lambda^n_t)\|
$$
then $ f_0(t) = f_n(0) = \|(\mu_0, \lambda_0)\| \leq 1 $ and by the
estimate \eqref{C1}
$$
f_{n + 1} (t) \leq 1 + C \int^t_0 f_n (s)^2 \, ds.
$$
Hence 
$$
f_n(t) \leq (1 - Ct)^{-1}, \quad t \leq C^{-1}
$$
for all $n$, so, setting $ T = (2C)^{-1} $, we have
\begin{equation}
\eqlabel{T}
\|(\mu^n_t, \lambda^n_t) \| \leq 2, \quad t \leq T.
\end{equation}
\page{20}
Next, set $ g_0(t) = f_0(t) $ and for $ n \geq 1 $
$$
g_n(t) = \|(\mu^n_t, \lambda^n_t) - (\mu^{n - 1}_t, \lambda^{n -
1}_t)\|.
$$
By the estimates \eqref{C2} and \eqref{T}, there is a
constant $ C < \infty $, depending only on $ \varphi $ and $B$, such
that
$$
g_{n + 1} (t) \leq C\int^t_0 g_n(s) \, ds, \quad t \leq T.
$$
Hence, by the usual arguments, $ (\mu^n_t, \lambda^n_t) $ converges in
$ \cal M_B \times \R $, uniformly in $ t \leq T $, to the desired local
solution, which is also unique. Moreover, for some constant $ C <
\infty $, depending only on $ \varphi $ and $ B $, we have
$$
\|(\mu_t, \lambda_t)\| \leq C, \quad t \leq T.
$$
\page{21}

It remains to show that $ \mu_t \geq  0 $ for all $t$. For this we need
the following result.

\begin{proposition}
\label{2.3}
Let
$$
(t,x) \mapsto f_t(x): [0,T] \times B \rightarrow \R
$$
be a bounded measurable function, having a bounded partial derivative
$ \partial f/\partial t $. Then for all $ t \leq T $ 
$$
\frac{d}{dt} \,\< f_t, \mu_t \> = \< \partial f/\partial t,
\mu_t \> + \< (f_t, 0), L^B (\mu_t, \lambda_t) \>.
$$
\end{proposition}
\begin{proof}
Fix $ t \leq T $ and set $ \lfloor s\rfloor_n = (n/t)^{-1} \lfloor
ns/t \rfloor $ and $ \lceil s \rceil_n = (n/t)^{-1} \lceil ns/t \rceil
$. Then
$$
\< f_t, \mu_t \> = \< f_0, \mu_0 \> + \int^t_0 \< \partial f/ \partial
s, \mu_{\lfloor s \rfloor_n} \> \, ds + \int^t_0 \< ( f_{\lceil s
\rceil_n}, 0), L^B (\mu_s, \lambda_s) \> \, ds
$$
and the proposition follows on letting $ n \rightarrow \infty $.
\end{proof}
\page{22}

For $ t \leq T $, set
$$
\theta_t(x) = \exp \int^t_0 \bigg( \int_E K (x,y) \mu_s (dy) + \lambda_s
\varphi(x) \bigg) \, ds
$$
and define $ G_t: \cal M_B \rightarrow \cal M_B $ by
$$
\< f, G_t (\mu)\> = \frac{1}{2} \int_{E \times E} (f \theta_t)( x + y) 1_{x +
y \in B} K (x,y) \theta_t (x)^{-1} \theta_t(y)^{-1} \mu (dx) \mu
(dy).
$$
We note that $ G_t (\mu) \geq 0 $ whenever $ \mu \geq 0 $ and, for
some $ C < \infty $, depending only on $ \varphi $ and $ B $, we have
$$
\|G_t(\mu) \| \leq C \|\mu\|^2, \,\|G_t(\mu) - G_t(\mu ')\| \leq C \|
\mu - \mu ' \| \,( \|\mu \| + \|\mu ' \|).
$$
Set $ \tilde \mu_t = \theta_t \mu_t $. By Proposition \ref{2.3}, for
all bounded measurable functions $f$, we have
$$
\frac{d}{dt} \< f, \tilde \mu_t \> = \< f \partial \theta/{\partial t} ,
\mu_t \> + \< (f \theta_t, 0), L^B (\mu_t, \lambda_t) \> = \< f,
G_t(\tilde \mu_t) \> .
$$
Define inductively a new sequence of measures $ \tilde \mu^n_t $ by
setting $ \tilde \mu^0_t = \mu_0 $ and, for $ n \geq 0 $
$$
\tilde \mu^{n + 1}_t = \mu_0 + \int^t_0 G_s(\tilde \mu^n_s) \, ds.
$$
\page{23}
By an argument similar to that used for the original iterative scheme,
we can show, first, and possibly for a smaller value of $ T >  0 $,
but with the same dependence, that $ \| \tilde \mu^n_t\| $ is bounded,
uniformly in $n$, for $ t \leq T $, and then that $ \|\tilde \mu^n_t -
\tilde \mu_t \| \rightarrow 0 $ as $ n \rightarrow \infty $. Since $
\tilde \mu^n_t \geq 0 $ for all $n$, we deduce $ \tilde \mu_t \geq  0 $
and hence $ \mu_t \geq 0 $ for all $ t \leq T $. This completes the
proof of Proposition \ref{2.2}
\end{proof}
\page{24}
We remark that the arguments used to prove Proposition \ref{2.2} apply
with no essential change to the case where the coagulation kernel is
time-dependent provided that \eqref{varphiK} holds uniformly in time.
We remark also that, in the iterative scheme
$$
\mu^0_t = \mu_0,\hspace*{.75in} \mu^{n + 1}_t \ll \mu_0 + \int^t_0 ( \mu^n_s + \mu^n_s * \mu^n_s ) \,
ds
$$
for all $ n \geq 0 $, so by induction we have
$$
\mu^n_t \ll \gamma_0 = \sum^\infty_{k = 1} \mu^{* k}_0
$$
where $ \mu^{* n}_0 $ denotes the $n$-fold convolution of 
$ \mu_0$. On letting $ n \rightarrow \infty $, we see
that, if $ (\mu_t, \lambda_t)_{t \geq 0} $
is the unique solution to \eqref{B}, then $ \mu_t
\ll \gamma_0 $. These remarks will be used in the proof of Proposition
\ref{4.2}.

We now fix $ \mu_0 \in \cal M^+ $ with $ \mu_0 \geq 0 $ and $ \< \varphi,
\mu_0 \> < \infty $. For each compact set $ B \subseteq E $, let
$$
\mu^B_0 = 1_{B} \mu_0, \quad \lambda^B_0 = \int_{E\backslash B}
\varphi(x) \mu_0 \, (dx)
$$
and denote by $ (\mu^B_t, \lambda^B_t)_{t \geq 0}$ the unique
solution to \eqref{B}, starting from 
$(\mu^B_0, \lambda^B_0)$, provided by Proposition \ref{2.2}.
We shall show in Proposition \ref{2.4} that for $ B \subseteq B ' $ we
have
$$
\mu^B_t \leq \mu^{B'}_t , \quad \< \varphi , \mu^B_t \> + \lambda^B_t
\geq \< \varphi, \mu^{B'}_t \> + \lambda^{B'}_t.
$$
We shall also show in Proposition \ref{2.5} that, for any local
solution $ (\nu_t)_{t < T} $ of the coagulation equation \eqref{C},
for all $t<T$
$$
\mu^B_t \leq \nu_t, \quad \< \varphi, \mu^B_t \> + \lambda^B_t \geq \<
\varphi, \nu_t \>.
$$
We now show how these facts lead to the proof of Theorem \ref{2.1}.

Set $ \mu_t = \lim_{B \uparrow E} \mu^B_t $ and $ \lambda_t
= \lim_{B \uparrow E} \lambda^B_t $. Note that
$$
\< \varphi, \mu_t \> = \lim_{B \uparrow E} \< \varphi, \mu^B_t \> \leq
\< \varphi , \mu_0 \> < \infty.
$$
So, by dominated convergence, using \eqref{varphiK}, for all bounded
measurable functions $f$,
$$
\int_{E \times E} f (x + y) 1_{ x + y \not\in B} K (x,y) \mu^B_t \,
(dx) \mu^B_t (dy) \rightarrow 0
$$
and we can pass to the limit in \eqref{B} to obtain
\page{26}
$$
\frac{d}{dt} \< f, \mu_t \> = \frac{1}{2} \int_{E \times E} \{ f(x + y) - f(x) - f(y)
\} K (x,y) \mu_t \, (dx) \mu_t (dy) - \lambda_t \< f \varphi, \mu_t \>.
$$
For any local solution $ (\nu_t)_{t < T} $, for all $ t < T $, 
$$
\mu_t \leq \nu_t, \quad \< \varphi, \mu_t \> + \lambda_t \geq \<
\varphi, \nu_t \>. 
$$
Hence, if $ \lambda_t = 0 $ for all $ t < T $, then
$ (\mu_t)_{t < T} $ is a local solution and moreover is the only local
solution on $ [0,T) $. If $ (\nu_t)_{t < T} $ is a strong local
solution, then
$$
\int^t_0 \< \varphi^2, \mu_s \> \, ds \leq \int^t_0 \< \varphi^2, \nu_s
\> \, ds < \infty
$$
for all $t<T$; this allows us to pass to the limit in \eqref{B} to obtain
\begin{equation}
\eqlabel{lambda}
\frac{d}{dt} \lambda_t = \lambda_t \< \varphi^2, \mu_t \>
\end{equation}
and to deduce from this equation that $ \lambda_t = 0 $ for all
$ t < T $; it follows that $ (\nu_t)_{t < T} $ is the only local
solution on $ [0, T) $.
\page{27}
Note that, for any local solution $ (\nu_t)_{t < T} $ 
\begin{multline*}
\int_E x 1_{x \leq n} \nu_t \, (dx) =
 \int_E x 1_{x \leq n} \nu_0 \, (dx) \\
+ \frac{1}{2} \int^t_0 \int_{E \times E} \{(x + y) 1_{x + y \leq n} - x 1_{x
\leq n} - y 1_{y \leq n} \} K (x,y) \nu_s (dx) \nu_s (dy) \,ds.
\end{multline*}
Hence, if $ (\nu_t)_{t < T} $ is strong and $ \varphi(x) \geq \varepsilon
x $ for all $x$, for some $ \varepsilon >  0 $, then by dominated
convergence the second term on the right tends to $ 0 $ as $ n
\rightarrow \infty $, showing that $ (\nu_t)_{t < T} $ is conservative.
\page{28}

Suppose now that $ \< \varphi^2, \mu_0 \> < \infty $ and set $ T = \<
\varphi^2 , \mu_0 \>^{-1} $. For any compact set $ B \subseteq E $ we
have
$$
\frac{d}{dt} \< \varphi^2, \mu^B_t \> \leq 
\frac{1}{2} \int_{E \times E} \{ \varphi (x
+ y)^2 - \varphi(x)^2 - \varphi(y)^2 \} K (x,y) \mu^B_t (dx) \mu^B_t
(dy) \leq \< \varphi^2, \mu^B_t \>^2
$$
so, for $ t  < T $,
$$
\< \varphi^2, \mu_t \> \leq \lim_{B \uparrow E} \< \varphi^2, \mu^B_t
\> \leq (T - t)^{-1}.
$$
Hence \eqref{lambda} holds and forces $ \lambda_t = 0 $ for $ t < T $ 
as above, so $ (\mu_t)_{t < T} $ is a strong local solution.

If $ \varphi^2 $ is sublinear, then
$$
\< \varphi^2, \mu_t \> \leq \< \varphi^2, \mu_0 \> < \infty.
$$
If, on the other hand, $ K (x,y) \leq \varphi(x) + \varphi(y) $, then 
$$
\eqalign{
\frac{d}{dt} \< \varphi^2, \mu^B_t \> &\leq \int_{E \times E} \varphi(x)
\varphi(y) (\varphi(x) + \varphi(y)) \mu^B_t (dx) \mu^B_t(dy)\cr
& \leq 2 \< \varphi, \mu^B_t \> \< \varphi^2, \mu^B_t \> 
 \leq 2 \< \varphi, \mu_0 \> \< \varphi^2, \mu^B_t \>,
}
$$
\page{29}
so
$$
\< \varphi^2, \mu_t \> \leq \exp \{ 2 \< \varphi, \mu_0 \> t \}.
$$
In either case we can deduce that $ (\mu_t)_{t \geq 0} $ is a strong
solution.
%\end{proof}

\page{30}
\begin{proposition}
\label{2.4}
Suppose $ B \subseteq B' $. Then, for all $ t \geq 0 $,
$$
\mu^B_t \leq \mu^{B'}_t, \quad \< \varphi , \mu^B_t \> + \lambda^B_t
\geq \< \varphi, \mu^{B'}_t \> + \lambda^{B'}_t.
$$
\end{proposition}
\page{31}

\begin{proof}
Set
$$
\eqalign{
\theta_t(x) & = \exp \int^t_0 \left(\int_E K (x,y)  \mu^B_s (dy) + \lambda^B_s
\varphi(x) \right) \, ds,\cr
\pi_t & = \theta_t (\mu^{B'}_t - \mu^B_t ), \cr
\chi_t &= \< \varphi, \mu^B_t \> + \lambda^B_t - \<  \varphi, \mu^{B'}_t
\> - \lambda^{B'}_t.
}
$$
Note that $ \pi_0 \geq 0 $ and $ \chi_0 = 0 $. By Proposition \ref{2.3},
for any bounded measurable function $f$
$$
\eqalign{
\frac{d}{dt} \< f, \pi_t \> &= \< f\, \partial \theta/\partial t, \mu^{B'}_t
- \mu^B_t \> + \< (f \theta_t, 0), L^{B'} (\mu^{B'}_t, \lambda^{B'}_t)
+ L^B (\mu^B_t, \lambda^B_t) \>\cr
&= \frac{1}{2} \int_{E \times E} f \theta_t (x + y) K (x,y) \cr
&\qquad\qquad\times(1_{x + y \in B'}
\mu^{B'}_t (dx) \mu^{B'}_t (dy) 
- 1_{x + y \in B} \mu^B_t (dx) \mu^B_t(dy)) \cr
&\quad+ \int_{E \times E} f \theta_t(x) (\varphi(x) \varphi(y) - K
(x,y)) \mu^{B'}_t (dx)( \mu^{B'}_t (dy) -  \mu^B_t(dy)) \cr
&\quad + \chi_t \int_E f \theta_t(x) \varphi(x) \mu^{B'}_t (dx).
}
$$
\page{32}
Also
$$
\frac{d}{dt} \, \chi_t = \frac{1}{2} \int_{E \times E} \{ \varphi(x) + \varphi(y) -
\varphi (x+y) \} K (x,y) (\mu^{B'}_t (dx) \mu^{B'}_t (dy) - \mu^B_t
(dx) \mu^B_t (dy)).
$$
So $ ( \pi_t, \chi_t) $ satisfies an equation of the form
$$
\frac{d}{dt} (\pi_t, \chi_t) = H_t (\pi_t, \chi_t) + (\alpha_t, 0)
$$
where $ H_t: \cal M_{B'} \times \R \rightarrow \cal M_{B'} \times \R $ is
linear, $ H_t (\pi, \chi) \geq 0 $ whenever $ (\pi, \chi) \geq 0 $, where $
\alpha_t \in \cal M_{B'} $ with $ \alpha_t \geq 0 $, and where we have
estimates, for $ t \leq 1 $,
$$
\| H_t (\pi, \chi)\|\leq C \|(\pi, \chi)\| \qquad\hbox{and}\qquad
\|\alpha_t \|  \leq C
$$
for some constant $ C < \infty $ depending only on $ \varphi $ and $
B' $.
\page{33}
Therefore we can apply the same sort of argument that we used for
non-negativity to see that $ \pi_t \geq 0 $ and $ \chi_t \geq 0 $ for all
$ t \leq 1 $, and then for all $ t < \infty $, as required.
Explicitly, $ H_t $ is given by
$$
\eqalign{\<(f,a),  H_t (\pi, \chi) \> &= 
 \frac{1}{2} \int_{E \times E} f \theta_t ( x + y) 1_{x + y \in B} K (x,y) 
\cr
&\qquad\qquad\times
(\theta_t (x)^{-1} \pi(dx) \mu^{B'}_t (dy)
+ \theta_t (y)^{-1} \mu^B_t (dx) \pi (dy)) \cr
&\quad + \int_{E \times E} f \theta_t(x) (\varphi(x) \varphi (y) - K
(x,y)) \mu^{B'}_t (dx) \theta_t(y)^{-1} \pi(dy)\cr
&\quad + \chi \int_E f \theta_t (x) \varphi(x) \mu^{B'}_t (dx) \cr
&\quad + \frac{1}{2} a \int_{E \times E} \{ \varphi (x) + \varphi(y) - \varphi(x + y) \} K
(x,y) \cr
&\qquad\qquad(\theta_t(x)^{-1} \pi(dx) \mu^{B'}_t (dy) 
+ \theta_t(y)^{-1} \mu^B_t (dx) \pi(dy))
}
$$
and $ \alpha_t $ is given by
$$
\< f, \alpha_t \> = \frac{1}{2} \int_{E \times E} f \theta_t (x + y) 1_{x + y
\in B'\backslash B} K (x,y) \mu^{B'}_t (dx) \mu^{B'}_t (dy).
%\qed
$$
\end{proof}
\page{34}
\begin{proposition}
\label{2.5}
Suppose that $ (\nu_t)_{t < T} $ is  a local solution of the coagulation
equation {\rm \eqref{C}}, starting from $ \mu_0 $. Then, for all compact
sets $ B \subseteq E $
and all $ t < T $,
$$
\mu^B_t \leq \nu_t, \quad \< \varphi, \mu^B_t \> + \lambda^B_t \geq \<
\varphi, \nu_t \>.
$$
\end{proposition}
\page{35}
\begin{proof}
Set
$$
\eqalign{
\theta_t(x) & = \exp \int^t_0 \left( \int_E K (x,y) \mu^B_s (dy) +
\lambda^B_s \varphi(x) \right) \, ds,\cr
\nu^B_t & = 1_B \nu_t,\cr
\pi_t & = \theta_t (\nu^B_t - \mu^B_t),\cr
\chi_t & = \< \varphi, \mu^B_t \> + \lambda^B_t - \<  \varphi, \nu_t \>.
}
$$
We have to show that $ \pi_t \geq 0 $ and $ \chi_t \geq 0 $. By an
obvious modification of Proposition \ref{2.3} we have, for all bounded
measurable functions $f$,
$$
\frac{d}{dt} \< f, \pi_t \> = \< f \, \partial \theta/\partial t, \nu^B_t -
\mu^B_t \> + \< f \theta_t 1_B, L(\nu_t)\> -\<(f\theta_t,0), L^B (\mu^B_t, \lambda^B_t)
\>.
$$
By \eqref{F}, we have
\page{36}
$$
\< \varphi, \nu_t \> \leq \< \varphi, \nu_0 \> + \frac{1}{2} \int^t_0 \int_{E
\times E} \{ \varphi(x + y) - \varphi(x) - \varphi(y) \} K (x,y) \nu_s
(dx) \nu_s (dy) \, ds.
$$
On the other hand
$$
\frac{d}{dt} ( \< \varphi, \mu^B_t \> + \lambda^B_t) = \frac{1}{2} \int_{E
\times E} \{ \varphi(x + y) - \varphi (x) - \varphi(y) \} K (x,y)
\mu^B_t(dx) \mu^B_t (dy).
$$
So $ \chi_t \geq \rho_t $ where
$$
\rho_t = \frac{1}{2} \int^t_0 \int_{E \times E} \{ \varphi (x) + \varphi (y) -
\varphi (x + y) \} K (x,y) ( \nu_s (dx) \nu_s(dy) - \mu^B_s (dx) \mu^B_s
(dy)) \, ds.
$$
Now $ (\pi_t, \rho_t) \in \cal M_B \times \R $ obeys a differential
equation of the form
$$ 
\frac{d}{dt}(\pi_t, \rho_t ) = H_t (\pi_t, \rho_t) + (\alpha_t, \beta_t)  
$$
where $ H_t: \cal M_B \times \R \rightarrow \cal M_B \times \R $ is linear, $
H_t (\pi, \rho) \geq 0 $ whenever $ (\pi, \rho) \geq 0 $, where $ \alpha_t
\geq 0 $, $ \beta_t \geq 0$ 
and we have estimates of the form
$$
\|H_t (\pi, \rho)\| \leq C \|(\pi, \rho)\|, \|\alpha_t\| \leq C,
|\beta_t| \leq C.
$$
\page{37}
It follows that $ \pi_t \geq 0 $ and $ \rho_t \geq 0 $, so also $
\chi_t \geq 0 $ as required.
Explicitly, $ H_t $ is given by
$$
\eqalign{
\< (f,a), H_t (\pi, \rho) \> &= 
 \frac{1}{2} \int_{E \times E} f \theta_t (x + y) 1_{ x + y \in B} K (x,y)
\cr
&\qquad\qquad\times(\theta_t (x)^{-1}\pi(dx) \nu^B_t(dy) 
+ \theta_t (y)^{-1} \mu^B_t (dx) \pi(dy))\cr
&\quad + \int_{E \times E} f \theta_t(x) (\varphi(x) \varphi(y) - K (x,y))
\theta_t (y)^{-1} \nu^B_t (dx) \pi(dy)\cr
&\quad + \rho \int_E f \theta_t(x) \varphi(x) \nu^B_t (dx)\cr
&\quad + \frac{1}{2} a \int_{E \times E} \{ 
\varphi(x) + \varphi(y) - \varphi(x +
y)\} K (x,y) \cr
&\qquad\qquad\times (\theta_t (x)^{-1} \pi(dx) \nu^B_t(dy) 
+ \theta_t (y)^{-1} \mu^B_t (dx) \pi (dy) )
}
$$
and $ \alpha_t , \beta_t $ are given by
$$
\eqalignbot{
\< f, \alpha_t \> &= \int_{E \times E} f \theta_t(x) (\varphi(x)
\varphi(y) - K (x,y)) 1_{y\not\in B} \nu^B_t(dx) \nu_t(dy) \cr
&\quad + (\chi_t - \rho_t) \int_E f \theta_t (x) \varphi(x) \nu^B_t(dx) \cr
&\quad  + \frac{1}{2} \int_{E \times E} f \theta_t ( x + y) 1_{ x + y \in B} K
(x,y) 1_{(x,y) \not\in B \times B} \nu_t (dx) \nu_t(dy),\cr
\beta_t & = \frac{1}{2} \int_{E \times E} \{ \varphi(x) + \varphi(y) - \varphi(x +
y) \} K (x,y) 1_{(x,y) \not\in B \times B} \nu_t (dx) \nu_t (dy).
}
\eqno
\qed
$$
\let\qed\relax
\end{proof}
\page{38}
This concludes the proof of Theorem \ref{2.1}.
\marginpar{added this line}
%\end{proof}
\page{39}

\section{An example of non-uniqueness}
\label{3}
We construct in this section an example of Smoluchowski's coagulation
equation having at least two solutions, both of which are moreover
conservative. 
\page{40}

Consider the system of differential equations
\begin{equation}
\eqlabel{m}
\frac{d}{dt} m_n(t) = - \lambda_n m_n (t) m_{n + 1} (t), \quad n = 1,2,
\ldots 
\end{equation}
For any solution we have
$$
m_n(t) = m_n(0) \exp\{- \lambda_n \int^t_0 m_{n + 1} (s) \, ds\}.
$$
Assume that $ m_n(0) \geq 0 $ for all $n$, then $ m_n(t) \geq 0 $ for
all $n$ and $t$. Note that, if $ m_N $ is given, and we consider the
system restricted to $ n \leq N - 1 $, then $ m_n $ is decreasing in $
m_N $ when $ N - n $ is odd, and increasing in $ m_N $ when $ N - n $
is even.
\page{41}

Fix $N$ and consider the case
$$
m_n(0) = \cases{
2^{-n}, & $n = 1, \ldots, 2N$\cr
0, & $ n = 2N + 1, \ldots $
}
$$
and 
$$
\lambda_n = 8^n \quad \text{ for all } n.
$$

\begin{proposition}
\label{3.1}
We have for all $n$ and $t$
$$
\eqalign{
m_{2n} (t) & \geq  \frac{1}{2} m_{2n} (0),\cr
m_{2n + 1}(t) & \leq m_{2n + 1} (0) \exp\{-4^{2n}t\}.
}
$$
\end{proposition}

\begin{proof}
Certainly $ m_{2N} (t) = m_{2N}(0) \geq \frac{1}{2} m_{2N}(0)$. Suppose that 
$n < N $ and 
$$ 
m_{2n + 2} (t) \geq \frac{1}{2} m_{2n + 2}(0)\quad \text{ for all } t.
$$ 
Then
$$
\eqalign{
m_{2n + 1}(t) & \leq m_{2n + 1} (0) \exp\{-\lambda_{2n + 1} m_{2n + 2}
(0) t/2\} \cr
& = m_{2n + 1}(0) \exp\{- 4^{2n}t\}
}
$$
so
$$
\int^\infty_0 m_{2n + 1} (t) \, dt \leq m_{2n + 1} (0) 4^{-2n} = \frac{1}{2}\,
8^{-2n}
$$
so
$$
m_{2n}(t) \geq m_{2n} (0) \exp\{-8^{-2n}\lambda_{2n}/2\} \geq \frac{1}{2}
m_{2n}(0).
$$
Hence the proposition follows by induction.
\end{proof}
\page{42}

We denote the solution just considered by $ m^{2N} $. The same
arguments establish corresponding inequalities for $ m^{2N + 1} $,
where the roles of even and odd are swapped. Now we let $ N
\rightarrow \infty $. For $ N \geq n $, $ m^{2N}_{2n}(t) $ is
decreasing in $N$ and $ m^{2N}_{2n + 1}(t) $ is increasing in $N$ for
all $n$. We set
$$
m^+_n (t) = \lim_{N \rightarrow \infty} m^{2N}_n(t).
$$
The integral equation
$$
m^{2N}_n(t) = m^{2N}_n (0) - \lambda_n \int^t_0 m^{2N}_n (s) m^{2N}_{n
+ 1} (s) \, ds
$$
passes to the limit to give
$$
m^+_n(t) = m^+_n (0) - \lambda_n \int^t_0 m^+_n(s) m^+_{n+1}(s)\, ds
$$
so $ m^+ $ is differentiable with
$$
\frac{d}{dt} m^+_n(t) = - \lambda_n m^+_n (t) m^+_{n + 1} (t).
$$
So $ (m^+_n : n \geq 1)$ solves the original system of equations. The
same argument produces another solution $ (m^-_n : n \geq 1) $ given
by
\page{43}
$$
m^-_n(t) = \lim_{N \rightarrow \infty} m^{2N + 1}_n (t).
$$
We have $ m^+_n (0) = m^-_n(0) = 2^{-n} $ for all $n$.
But
$$
\eqalign{
m^+_{2n}(t) & \geq \frac{1}{2} m^+_{2n} (0), \cr
m^-_{2n} (t) &\leq \exp\{-4^{2n - 1}t\} m^-_{2n} (0)
}
$$
for all $n$ and $t$, so $ m^+ \neq m^- $.
\page{44}

We now use the solutions $ m^+ $ and $ m^-$ to construct an example of 
Smoluchowski's coagulation equation having at least two conservative
solutions. Let $ x_1, x_2, \ldots $ be an increasing sequence in $ (0,
\infty) $ which is linearly independent over $ \Z $. For 
$$
x = x_n + (k_1 x_1 + \ldots + k_{n - 1} x_{n - 1}), \quad k_1, \ldots,
k_{n - 1} \in \Z
$$
we write $ n(x) = n $. Denote by $I$ the set of all such $x$ and
define $n(x)=0$ if $x\notin I$. Define $ K: E \times E \rightarrow [0, \infty) $ by
$$
K (x,y) = \cases{
\lambda_n & if $\{ n(x),n(y)\} =\{n, n + 1\} $ and $n\geq 1,$ \cr
0 & otherwise.
}
$$
Set
$$
\mu_0 = \sum^\infty_{n = 1} \varepsilon_{x_n} 2^{-n}
$$
and consider Smoluchowski's coagulation equation
$$
\frac{d}{dt} \< f, \mu_t \> = \frac{1}{2} \int_{E \times E} \{ f(x + y) - f(x)
- f(y) \} K (x,y) \mu_t (dy) \mu_t (dx)
$$
starting from $ \mu_0 $.
\page{45}

According to the definition made in \S2, for a solution, we require
$$
\int^t_0 \int_{B \times E} K (x,y) \mu_s (dx) \mu_s (dy)\, ds < \infty
$$
for all $t$ and all compact sets $ B \subseteq E $, and, for all
bounded measurable functions $f$ of compact support
$$
\<f, \mu_t \> = \< f, \mu_0 \> +
\frac{1}{2} \int^t_0 \int_{E \times E} \{ f(x + y) - f(x) - f(y) \} K (x,y)
\mu_s (dx) \mu_s (dy) \, ds.
$$
Consider for $ n = 1, 2, \ldots, $
$$
m_n(t) = \mu_t ( \{ x: n(x) = n \} ).
$$
Take $ f(x) = 1_{n(x) = n, x \leq k} $ and let $ k \rightarrow \infty $
to obtain
$$
\frac{d}{dt} m_n(t) = - \lambda_n m_n (t) m_{n + 1} (t).
$$
We deduce that any solution $ (\mu_t)_{t \geq 0} $ of the coagulation
equation gives rise to a solution $ (m_n(t) : n \geq 1)_{t \geq 0} $
of the system \eqref{m}.
\page{46}
On the other hand, for any solution $ (m_n(t): n \geq 1)_{t \geq 0} $ of this
system, we obtain a solution $ (\mu_t)_{t \geq 0} $ of the coagulation
equation by
$$
\frac{d}{dt} \mu_t (\{ x \} ) = - 
(\lambda_{n-1} m_{n - 1} (t)+ 
\lambda_n m_{n + 1} (t)) 
\mu_t ( \{ x \} )  
+ \frac{1}{2} \lambda_{n - 1} \sum_{\scriptstyle y, z \in I \atop y + z
= x} \mu_t (\{ y \} ) \mu_t ( \{ z \} )
$$
whenever $ n(x) = n $. Hence the coagulation equation has two distinct
solutions $ (\mu^+_t)_{t \geq 0} $ and $ (\mu^-_t)_{t \geq 0} $
corresponding to $ (m^+_t)_{t \geq 0} $ and $ (m^-_t)_{t \geq 0}$. We
now show these solutions are conservative.
The idea of the proof is to show that the proportion of original particles making at
least $2n$ jumps falls off geometrically in $n$.
\page{47}
\begin{proposition}
\label{3.2}
Suppose that $ \mu_0 $ has finite total mass density. Then the
solution $ (\mu^+_t)_{t \geq 0} $ is conservative.
\end{proposition}
\begin{proof}
For $ x = k_1 x_1 + \ldots + k_{n - 1} x_{n - 1} + x_n \in I $,
define 
$$
k_m (x) = \cases{
k_m & if $ m < n $ \cr
1 & if $ m = n $ \cr
0 & if $ m >  n $. \cr
}
$$
Note that $ k_m: I \rightarrow \Z^+ $ is additive. For $ m \leq n $
consider 
$$
v_{m,n} (t) = \int_I k_m(x) 1_{n(x) \leq n} \mu_t (dx).
$$
Then $v_{m,n}(t)$ is non-increasing in $t$.
Note that
$$
\int_{I \times I} k_m (x) 1_{n(x) \leq n} K (x,y) \mu_t (dx) \mu_t
(dy) \leq 8^{n}v_{m,n}(t) <\infty 
$$
\marginpar{ms. illegible at ??}
so, by dominated convergence, since $ (\mu^+_t)_{t \geq 0} $ is a solution,
$$
\frac{d}{dt} v_{m,n} (t) = - \lambda_n m_{n + 1} (t) \int_I
k_m(x) 1_{n(x) = n} \mu_t (dx).
$$
Note also
$$
\frac{d}{dt} \int_I k_m(x) 1_{n(x) = n} \mu_t (dx) \leq -\frac{d}{dt}
v_{m,n - 1}(t)
$$
so
$$
\int_I k_m(x) 1_{n(x) = n} \mu_t (dx) \leq r_{m,n - 1}(t)
$$
where $ r_{m, n}(t) = 2^{-m} - v_{m,n}(t) $. Now for $n$ even
\page{48}
$$ 
r_{m,n}(t) \leq r_{m, n - 1}(t) \int^t_0 \lambda_n m_{n + 1} (s) \,
ds \leq \frac{1}{2} r_{m, n - 1}(t)
$$ 
so $ r_{m,n}(t) \rightarrow 0 $ as $ n \rightarrow \infty $. 
Hence $ v_{m,n}(t) \rightarrow 2^{-m} $ as $ n \rightarrow
\infty $ and hence
$$
\int_I x \mu_t(dx) = \sum_m x_m \int_I k_m(x) \mu_t(dx) = \sum_m x_m
2^{-m} = \int_I x \mu_0 (dx).
$$
\end{proof}

We make some remarks on the relation between this example and
the results of \S2. The construction of the example makes it
insensitive to the additive structure of $E$. We require very little
of the sequence $ (x_n)_{n \geq 1} $. By taking $ x_n \approx 8^n $ we
can satisfy the condition 
$$
K (x,y) \leq x
$$
but we get
$$
\int_E x\mu_0 (dx) = \infty.
$$
On the other hand, by taking $ x_n \approx \alpha^n $, for some $
\alpha < 2 $, we get
\page{49}
$$
\int_E x \mu_0 (dx) < \infty
$$
but the relation $ K (x,y) \leq C xy $ for all $ x, y \in E $ does not
hold for any $ C < \infty $. Thus, however we choose $ (x_n)_{n \geq
1} $, we cannot regard $ (\mu^{\pm}_t)_{t \geq  0} $ as a strong
solution, even in small time. This is
of course implied also by the uniqueness of strong solutions
established in \S2.
\page{50}

It would be nice to find an example of this type where the initial 
mass distibution is supported on $\N$. It may be that for integers
$x_n\rightarrow \infty$ sufficiently fast, the analogous equation
exhibits the same sort of behaviour. However we have not established
whether this is true.
\section{Hydrodynamic limit for the stochastic coalescent}
\label{4}
In this section we shall prove some limit theorems for the stochastic
coalescent. There are two main results. In Theorem \ref{4.1},
generalizing a result of Jeon \cite{Je}, we prove a tightness result
for the stochastic coalescent, which implies a general existence
theorem for solutions of Smoluchowski's equation. Then, in Theorem
\ref{4.4}, we prove weak convergence of the stochastic coalescent to
any strong solution of Smoluchowski's equation. The methods used are
mostly standard tools from the theory of weak convergence on Skorohod
spaces. The problem-specific idea which leads to Theorem \ref{4.4} is
the construction of a coupled family of particle systems, converging
to the stochastic coalescent, in direct analogy with the method of
\S2. A version of this idea was also discovered independently by 
Kurtz \cite{K}. The case of a discrete mass distribution may also be
treated using a differential equation approach instead of weak
convergence:~this is simpler and more effective, establishing
convergence at an exponential rate in the number of particles.
\page{51}
The particle system we consider has been considered, in various special
cases, by many others. In particular, it was considered in full
generality by Marcus \cite{M} and  Lushnikov \cite{Lu}. Recall that $
E = (0, \infty) $ and that the coagulation kernel $K$ is a symmetric
measurable function $ K: E \times E \rightarrow [0, \infty)$.

Let $ X_0 $ be a finite, integer-valued measure on $E$. We can write $
X_0 $ as a sum of unit masses
$$
X_0 = \sum^m_{i = 1} \varepsilon_{x_i}
$$
for some $ x_1, \ldots, x_m \in E $. We think of $ X_0 $ as
representing a system of $m$ particles, labelled by their masses $
x_1, \ldots, x_m $. A Markov process $ (X_t)_{t \geq 0} $ of finite,
integer-valued measures on $E$ can be constructed as follows:~for each
pair $ i < j $, take an independent exponential random time $ T_{ij} $
of parameter $ K(x_i, x_j) $ and set $ T = \min_{i < j} T_{ij}$; set $
X_{t} = X_0 $ for $ t < T $ and set
$$
X_T = X_0 - \varepsilon_{x_i} - \varepsilon_{x_j} + \varepsilon_{x_i + x_j}
\quad \text{ if } T = T_{ij};
$$
then begin the construction afresh from $ X_T $. In this process, each
pair of particles $ \{ x_i, x_j \} $ coalesces at rate $ K(x_i, x_j) $
to form a new particle $ x_i + x_j $. We call $ (X_t)_{t \geq 0} $ a
{\em  stochastic coalescent with coagulation kernel} $K$.
\page{53}

Denote by $d$ some metric on $\cal M$ which is compatible with the topology of
weak convergence, that is to say $ d(\mu_n, \mu) \rightarrow 0 $ if
and only if $ \< f, \mu_n \> \rightarrow \< f, \mu \> $ for all
bounded continuous functions $ f: E \rightarrow \R $. We choose $d$ so
that $ d (\mu, \mu ') \leq \| \mu - \mu ' \| $ for all $ \mu, \mu '
\in\cal M $.
\page{54}
When the class of functions $f$ is restricted to those of bounded
support we get a weaker topology, also metrizable, and we denote by $
d_0 $ some compatible metric, with $ d_0 \leq d $.

\page{53}
The following result is a first attempt at proving weak convergence
for the stochastic coalescent. It is less than satisfactory because it
does not enable us to show uniqueness of limits. We include it here,
partly as a warm-up for the more intricate arguments used later, and
partly because it provides the best result on global existence of
solutions to Smoluchowski's equation that we know. A version of the
result where $ \mu_0 $ is supported on $ \N $ and where $ \varphi(x) =
x $ has been proved already by Jeon \cite{Je}.

\page{55}
\begin{theorem}
\label{4.1}
Let $ K: E \times E \rightarrow [0, \infty) $ be a symmetric
continuous function and let $ \mu_0 $ be a measure on $E$. Assume
that, for some continuous sublinear function $ \varphi: E \rightarrow
(0, \infty) $, 
$$
K (x,y) \leq \varphi (x) \varphi(y), \quad \text{ for all } x, y \in E,
$$
$$
\varphi(x)^{-1} \varphi(y)^{-1} K (x,y) \rightarrow 0, \quad \text{ as
} (x,y) \rightarrow \infty.
$$
Assume also that $ \< \varphi, \mu_0 \> < \infty $. Let $ (X^n_t)_{t \geq 0} $ be a
sequence of stochastic coalescents, with coagulation kernel $K$. Set $
\tilde X^n_t = n^{-1} X^n_{n^{-1}t} $ and suppose that
$$
d_0 (\varphi \tilde X^n_0, \varphi \mu_0) \rightarrow 0
$$
as $ n \rightarrow \infty $ and that, for some constant $ \Lambda <
\infty $, for all n
$$
\< \varphi, \tilde X^n_0 \> \leq \Lambda.
$$
Then the sequence of laws of $ \tilde X^n $ on $ D([0, \infty), (\cal M,
d_0)) $ is tight. Moreover, for any weak limit point $X$, almost
surely, $ (X_t)_{t \geq 0} $ is a solution of Smoluchowski's coagulation
equation {\rm \eqref{C}}. In particular, this equation has at least one
solution.
\end{theorem}
\page{56}
\begin{proof}
For an integer-valued measure $ \mu $ on $E$, denote by $ \mu^{(1)} $
the integer-valued measure on $ E \times E $ given by
$$
\mu^{(1)} (A \times A') = \mu(A) \mu(A') - \mu(A \cap A').
$$
(This is simply the counting measure for ordered pairs of masses of
distinct particles.) Similarly, when $ n \mu $ is an integer-valued
measure, set
$$
\mu^{(n)} (A \times A') = \mu(A) \mu(A') - n^{-1} \mu (A \cap A')
$$
and note that $ n^2 \mu^{(n)} = (n \mu)^{(1)} $. For a bounded
measurable function $f$ on $E$, set
$$
\eqalign{
L^{(n)} (\mu)(f) & = \< f, L^{(n)} (\mu) \> \cr
& = \frac{1}{2} \int_{E \times E} \{ f(x + y) - f(x) - f(y) \} K (x,y)
\mu^{(n)} (dx, dy),\cr
Q^{(n)} (\mu)(f) & = \frac{1}{2} \int_{E \times E} \{ f(x + y) - f(x) - f(y)
\}^2 K (x,y) \mu^{(n)} (dx, dy).
}
$$
Then
$$
M^{f,n}_t = \< f, X^n_t \> - \< f, X^n_0 \> - \int^t_0 L^{(1)}
(X^n_s)(f) \, ds
$$
is a martingale, having previsible increasing process
$$
\< M^{f,n} \>_t = \int^t_0 Q^{(1)} (X^n_s)(f) \, ds.
$$
Set $ \tilde M^{f,n}_t = n^{-\frac{1}{2}} M^{f,n}_{n^{-1}t} $, then we have
\begin{equation}
\eqlabel{M}
\< f, \tilde  X^n_t \>  = \< f, \tilde X^n_0 \> + n^{-\frac{1}{2}} \tilde
M^{f,n}_t + \int^t_0 L^{(n)} (\tilde X^n_s) (f)\, ds, 
\end{equation}
$$
\< \tilde M^{f,n} \>_t  = \int^t_0 Q^{(n)} (\tilde X^n_s) (f) \, ds.
$$
Since $ \varphi $ is subadditive, we have $ \< \varphi, \tilde X^n_t
\> \leq \Lambda $ for all $n$ and $t$. Hence by \eqref{varphiK} 
$$
\eqalign{
|L^{(n)} (\tilde X^n_t) (f) | & \leq 2 \|f\| \Lambda^2,\cr
Q^{(n)} ( \tilde X^n_t)(f) & \leq 4 \|f\|^2 \Lambda^2.
}
$$
Assume that $|f|\leq\varphi\wedge 1$. Then
$$
|\<f,\tilde X^n_t\>|\leq\Lambda
$$
for all $t$, so we have compact containment.
Moreover, by Doob's $L^2-$inequality, for all $s<t$,
$$
\E\sup_{s\leq r\leq t} |\tilde M^{f,n}_r - \tilde M^{f,n}_s|^2
\leq4\E\int^t_s Q^{(n)}(\tilde X^n_r)(f)dr\leq 16\Lambda^2(t-s)
$$
so
$$
\E\sup_{s\leq r\leq t}|\<f,\tilde X^n_r-\tilde X^n_s\>|^2
\leq C\{(t-s)^2+n^{-1}(t-s)\}
$$
where $C<\infty$ depends only on $\Lambda$.
Hence, by a standard tightness criterion, the laws of the sequence
$\<f,\tilde X^n\>$ on $D([0,\infty),\R)$ are tight.
See for example \cite{EK}, Corollary 7.4.
We note the bound
$$
\|(\varphi\wedge 1)\tilde X^n_t\|\leq \<\varphi,\tilde X^n_t\>\leq\Lambda
$$
for all $t$.
Hence we can apply Jakubowski's criterion \cite{J} to see that the laws of the
sequence $ (\varphi \wedge 1) \tilde X^n $ on $ D([0, \infty), 
\cal M_{[0, \infty]}) $ are tight. 
By consideration of subsequences and a theorem
of Skorohod, see, for example, Pollard \cite{P}, Chapter IV, it
suffices from this point on to consider the case where $ (\varphi
\wedge 1) \tilde X^n $ converges almost surely in $ D([0, \infty),
\cal M_{[0, \infty]}) $, with limit $ (\varphi \wedge 1)X $ say. We denote
also by $X$ the process in $ \cal M_E $ obtained by restriction of
measures. Note that
$$
\| \tilde X^n_t - \tilde X^n_{t-} \| \leq 3/n
$$
so $ X \in C([0, \infty), \cal M_E) $. Moreover $ \varphi^\delta  = \varphi
1_{(0, \delta]} $ is subadditive, so
$$
\< \varphi^\delta , \tilde X^n_t \> \leq \< \varphi^\delta, \tilde
X^n_0 \> \leq \< \varphi^\delta, \mu_0 \> + | \< \varphi^\delta,
\tilde X^n_0 - \mu_0 \> |
$$
and so, given $ \varepsilon >  0 $, we can find $ \delta >  0 $ so that 
$$
\sup_n \, \sup_t \< \varphi^\delta, \tilde X^n_t \> < \varepsilon.
$$
Given a continuous bounded function $ f: E \rightarrow \R $ of bounded
support, we can write $ f = f_1 + f_2 $ where $ f_1 $ is continuous of
compact support and $ f_2 $ is supported in $ (0, \delta) $ with $
\|f_2 \| \leq \|f\| $. Then
\page{59}
$$
\limsup_{n \rightarrow \infty} \, \sup_{s \leq t} \<
\varphi f, \tilde X^n_s - X_s  \> 
 \leq \lim_{n \rightarrow \infty} \sup_{s \leq t} \<
\varphi f_1, \tilde X^n_s - X_s \>
 + \|f\| \sup_n \, \sup_s \< \varphi^\delta_, \tilde X^n_s \> 
 \leq \varepsilon \|f\|.
$$
Since $f$ and $ \varepsilon $ were arbitrary, this shows that
\begin{equation}
\eqlabel{dO}
\sup_{s \leq t} d_0 (\varphi  \tilde X^n_s, \varphi X_s) \rightarrow 0
\quad \text{ a.s. }
\end{equation}

We now wish to pass to the limit in \eqref{M}. Let us suppose for now
that $ f: E \rightarrow \R $ is continuous and of compact support $B$.
Then, as $ n \rightarrow \infty $
$$
\eqalign{ 
\E(\sup_{s \leq t} | n^{-\frac{1}{2}} \tilde M^{f,n}_s |^2 )
&\leq \frac{4}{n} \E \< \tilde M^{f,n}\>_t \leq \frac{16 \Lambda^2
\|f\|^2}{n} \rightarrow 0,\cr 
|(L - L^{(n)}) (\tilde X^n_s)(f)| &=
\frac{1}{n} \left|\int_E \{ f(2x) - 2f(x)\} K(x,x) \tilde X^n_s
(dx)\right| \cr 
& \leq\frac{3\|f\|}{n} \int_{B \cup 2B} \varphi(x)^2
\tilde X^n_s (dx) \cr 
& \leq \frac{3\|f\|}{n}\, \|\varphi 1_{B \cup
2B}\| \< \varphi, \tilde X^n_0 \> \rightarrow 0.  }
$$
Hence it will suffice to show that as $ n \rightarrow \infty $
\begin{equation}
\eqlabel{LX}
\sup_{s \leq t} |\< f, L(\tilde X^n_s) - L(X_s) \> | \rightarrow 0 \quad
\text{ a.s. }
\end{equation}
\page{60}
where we recall that
$$
\< f, L(\mu) \> = \frac{1}{2} \int \{ f(x + y) - f(x) - f(y) \} K (x,y)
\mu(dx) \mu(dy).
$$
Given $ \delta >  0 $ and $ N < \infty $, we can write $ K = K_1 + K_2
$, where $ K_1 $ is continuous of compact support and where $ 0 \leq
K_2 \leq K $ and $ K_2 $ is supported on
$$
F_1 \cup F_2 \cup F_2 = \{(x,y): x \leq \delta \} \cup \{ (x,y) : y
\leq \delta \} \cup \{ (x,y): | (x,y) | \geq N \}.
$$
Then, with an obvious notation, by \eqref{dO}
$$
\sup_{s \leq t} | \< f, L_1 (\tilde X^n_s) - L_1 (X_s) \> |
\rightarrow 0 \quad \text{ a.s.}
$$
whereas, for $ K_2 $, we use the estimates
$$
\eqalign{
\|K1_{F_1} \mu \otimes \mu \|  = \|K 1_{F_2} \mu \otimes \mu \|
& \leq \< \varphi, \mu \> \< \varphi^\delta, \mu \>,\cr
\|K1_{F_3} \mu \otimes \mu \| & \leq \beta_N \< \varphi, \mu \>^2,
}
$$
where $ \beta_N = \sup_{|(x,y) | \geq N} \varphi(x)^{-1} \varphi(y)^{-1} K
(x,y)$.
\page{61}
Now,
$$
\< \varphi, \tilde X^n_t \> \leq \< \varphi, \tilde X^n_0 \>, \quad \<
\varphi^\delta, \tilde X^n_t \> \leq \< \varphi^\delta, \tilde X^n_0
\>
$$
and, given $ \varepsilon >  0 $, we can find $ \delta $ and $N$ so that
$$
\eqalign{
\< \varphi^\delta, \tilde X^n_0 \> &\leq \tfrac{1}{3} \varepsilon \Lambda^{-1},
\quad \text{ for all } n,\cr
\< \varphi^\delta, \mu_0 \> &\leq \tfrac{1}{3} \varepsilon \Lambda^{-1}, \cr
\beta_N &\leq \tfrac{1}{3} \varepsilon \Lambda^{-2}.
}
$$
Then
$$
|\< f, L_2 (\tilde X^n_t ) \> | \leq \varepsilon , \quad |\< f, L_2
(X_t) \> | \leq \varepsilon
$$
for all $n$ and $t$. Hence
$$
\limsup_{n \rightarrow \infty} \, \sup_{s \leq t} | \< f, L( \tilde
X^n_s) - L(X_s) \> | \leq 2 \varepsilon.
$$
But $ \varepsilon $ was arbitrary, so \eqref{LX} is proved.
\page{62}
Hence we can let $ n \rightarrow \infty $ in  \eqref{m} to obtain
$$
\< f, X_t \> = \< f, X_0 \> + \int^t_0 \< f, L(X_s) \> \, ds
$$
for all continuous functions $ f: E \rightarrow \R $ of compact
support. By using the bounds \eqref{varphiK} and $ \< \varphi, X_t
\> \leq \Lambda $, and a straightforward limit argument, we can extend
this equation to all bounded measurable functions $f$. In particular,
almost surely, $X$ is a solution of Smoluchowski's equation, in the
sense of \S2.
\end{proof}
\page{63}

A corollary of Theorem of \ref{4.1} is that, whenever we know
Smoluchowski's equation has at most one solution, then, under the
hypotheses of Theorem \ref{4.1}, we can deduce, for all $t$
$$
\sup_{s \leq t} d_0 (\tilde X^n_s, \mu_s) \rightarrow 0
$$
in probability as $ n \rightarrow \infty $, for the solution $
(\mu_t)_{t \geq 0} $ provided by Theorem \ref{4.1}. However, we can
only prove uniqueness of solutions in the presence of a strong
solution. So we prefer to formulate our main limit result, Theorem
\ref{4.4}, in that context, when a new approach allows certain other
hypotheses to be relaxed.

For the remainder of this section we will assume that we have chosen a
continuous sublinear function $ \varphi: E \rightarrow (0, \infty) $
and that $K$ satisfies 
\marginpar{Labelled this `varphiK.2'}
\begin{equation}
\eqlabel{varphiK.2}
K (x,y) \leq \varphi(x) \varphi(y) \quad \text{ for all } x, y \in E.
\end{equation}
\page{66}

Our further analysis of the stochastic coalescent will rest on an
approximation by a coupled family of Markov processes $ (X^B_t,
\Lambda^B_t)_{t \geq 0} $, indexed by sets $ B \subseteq E $
which we now describe. Each process $ (X^B_t)_{t \geq 0} $ will take
values in the finite integer-valued measures on $E$, whereas $
(\Lambda^B_t)_{t \geq 0} $ will be a non-decreasing process in $ [0,
\infty) $. Let us suppose given initial values $ (X^B_0, \Lambda^B_0)
$, for all $B$, such that $ X^B_0 $ is supported in $B$ and such that,
whenever $ B \subseteq B' $
$$
X^B_0 \leq X^{B '}_0 , \quad \< \varphi, X^B_0 \> + \Lambda^B_0 \geq
\< \varphi, X^{B'}_0 \> + \Lambda^{B'}_0.
$$
Write $ X_0 = X^E_0 $ as a sum of unit masses
$$
X_0 = \sum^m_{i = 1} \varepsilon_{x_i}
$$
where $ x_1, \ldots, x_m \in E $. There is a unique increasing map
$$
B \mapsto I(B) \subseteq \{ 1, \ldots, m \} 
$$
such that
$$
X^B_0 = \sum_{i \in I(B)} \varepsilon_{x_i}.
$$
Set 
$$
\nu^B = \Lambda^B_0 - \sum_{j \not\in I(B)} \varphi (x_j).
$$
Note that $ \nu^B $ decreases as $B$ increases
\page{67}
and that $ \nu^E = \Lambda^E_0 \geq 0 $. For $ i < j $, take independent
exponential random variables $ T_{ij} $ of parameter $ K(x_i, x_j) $.
Set $ T_{ji} = T_{ij} $. Also, for $ i \neq j $, take independent
exponential random variables $ S_{ij} $ of parameter $ \varphi(x_i)
\varphi(x_j) - K(x_i, x_j) $. We can construct, independently
for each $i$, a family of independent exponential random variables $
S^B_i $, increasing in $B$, with $ S^B_i $ having parameter $ \varphi(x_i)
\nu^B $. Set
$$
T^B_i = \min_{j \not\in I(B)} (T_{ij} \wedge S_{ij}) \wedge S^B_i 
$$
and note that $ T^B_i $ is an exponential random variable of parameter
$$
\sum_{j \not\in I(B)} \varphi(x_i) \varphi(x_j) + \varphi(x_i) \nu^B =
\varphi(x_i) \Lambda^B_0.
$$
For each $B$, the random variables
$$
(T_{ij}, T^B_i: i, j \in I(B) , i < j )
$$
form an independent family, whereas, for $ i \in I(B) $ and $ j \not\in
I(B) $, we have
$$
T^B_i \leq T_{ij}
$$
\page{67}
and for $ B \subseteq B' $ and all $i$ we have
$$
T^B_i \leq T^{B'}_i.
$$
\marginpar{Added period}
Now set
$$
T = (\min_{i < j} T_{ij}) \wedge (\min_i T^{\emptyset}_i).
$$
We set $ (X^B_t, \Lambda^B_t ) = (X^B_0, \Lambda^B_0) $ for $ t < T $
and set
$$
(X^B_t, \Lambda^B_T )= \cases{
(X^B_0 - \varepsilon_{x_i} - \varepsilon_{x_j} + \varepsilon_{x_i + x_j},
\Lambda^B_0)& if $ T = T_{ij}, i, j \in I(B), x_i + x_j \in B $\cr
(X^B_0-\varepsilon_{x_i} - \varepsilon_{x_j}, \Lambda^B_0 + \varphi(x_i +
x_j))& if $ T = T_{ij}, i, j \in I(B) , x_i + x_j \not\in B $ \cr
(X^B_0 - \varepsilon_{x_i}, \Lambda^B_0 + \varphi(x_i))& if $ T = T^B_i,
i \in I(B) $\cr
(X^B_0, \Lambda^B_0) & otherwise.
}
$$
It is straightforward to check that $ X^B_T $ is supported on $B$ and,
for $ B \subseteq B' $
$$ 
X^B_T \leq X^{B'}_T, \quad \< \varphi , X^B_T \> + \Lambda^B_T \geq \<
\varphi, X^{B'}_T \> + \Lambda^{B'}_T.
$$
Hence we can repeat the above construction independently from time
$T$, to obtain a family of Markov processes $ (X^B_t, \Lambda^B_t)_{t
\geq 0} $ such that $ X^B_t $ is supported on $B$ and for $ B
\subseteq B' $ and all $t$ 
\begin{equation}
\eqlabel{XI}
X^B_t \leq X^{B'}_t, \quad\< \varphi, X^B_t \> + \Lambda^B_t \geq \<
\varphi, X^{B'}_t \> + \Lambda^{B'}_t .
\end{equation}
\page{68}
\noindent At the outset, we assumed that both $ X^B_0 $ and $ \Lambda^B_0 $ were
given, for all $B$. From now on we shall suppose simply that $ X_0 =
X^E_0 $ is given and take
$$
X^B_0 = 1_B X_0 , \quad \Lambda^B_0 = \< \varphi 1_{B^c}, X_0 \>.
$$
Of course these relations do not remain valid as time evolves.

For each fixed $B$, the process $ (X^B_t, \Lambda^B_t)_{t \geq 0}$ may
be regarded as a finite state-space Markov chain having three sorts of
transition. Each pair of particles $ x_i, x_j $ in $ X^B_t $ is, at
rate $ K(x_i, x_j) $, removed; if $ x_i + x_j \in B $, the merged
particle is added to $ X^B_t $, if not, $ \varphi(x_i + x_j) $ is
added to $ \Lambda^B_t $. Also, each particle $ x_i $ in $ X^B_t $ is,
at rate $ \varphi(x_i) \Lambda^B_t $, removed and $ \varphi(x_i) $
added to $ \Lambda^B_t $. In particular, for the choice of initial
values made above, $ \Lambda^E_t = 0 $ for all $t$ and $ X_t = X^E_t $
is simply the stochastic coalescent with coagulation kernel $K$ with
which we began.
\page{69}

We now proceed to identify some martingales associated with $ (X^B_t,
\Lambda^B_t)_{t \geq 0} $. Recall that, when $ n \mu $ is an
integer-valued measure on $E$, we denote by $ \mu^{(n)} $ the measure
on $ E \times E $ characterized by
$$
\mu^{(n)} (A \times A') = \mu(A) \mu(A') - n^{-1} \mu(A \cap A').
$$
Given \page{70} an integer-valued measure $ \mu $ on $E$ and given $ \lambda
\geq 0 $, define for any bounded measurable function $f$ on $E$ and
for $ a \in \R $
\begin{multline*}
L^{B,(1)} (\mu, \lambda) (f,a) = \< (f, a), L^{B,(1)} (\mu, \lambda)
\> = \\
\frac{1}{2} \int_{E \times E} \{ f(x + y) 1_{x + y \in B} + a \varphi(x +
y) 1_{ x + y \not\in B} - f(x) - f(y) \} K (x,y) \mu^{(1)} (dx, dy)\\
+ \lambda \int_E \{ a \varphi (x) - f (x) \} \varphi(x) \mu (dx)
\end{multline*}
and
\begin{multline*}
Q^{B, (1)} (\mu, \lambda)(f, a)
= \\
\frac{1}{2} \int_{E \times E} \{ f(x + y) 1_{x + y \in B} + a \varphi(x + y)
1_{ x + y \not\in B} -f(x) - f(y) \}^2 K (x,y) \mu^{(1)} (dx, dy)\\
+ \lambda \int_E \{ a \varphi (x) - f(x)\}^2 \varphi(x) \mu (dx).
\end{multline*}
\page{71}
Then, for all $f$ and $a$,
\begin{equation}
\eqlabel{L}
M_t = \< f, X^B_t \> + a \Lambda^B_t - \< f, X^B_0 \> - a \Lambda^B_0
- \int^t_0 L^{B(1)} (X^B_s, \Lambda^B_s)(f,a) \, ds
\end{equation}
is a martingale, having previsible increasing
process
\begin{equation}
\eqlabel{Q}
\< M \>_t = \int^t_0 Q^{B,(1)} (X^B_s, \Lambda^B_s ) (f,a) \, ds.
\end{equation}
Recall \page{72} from \S2 that, for $ B \subseteq E $ compact, we
denote by $ \cal M_B $ the space of finite signed measures supported on $B$ and
we define $ L^B: \cal M_B \times \R \rightarrow \cal M_B \times \R $ by the
requirement
\begin{multline*}
\<(f,a), L^B (\mu, \lambda) \> = \\
 \frac{1}{2} \int_{E \times E} \{ f(x + y) 1_{ x + y \in B} -  a \varphi( x
+ y) 1_{ x + y \not\in B} - f(x) - f(y) \} K (x,y) \mu (dx) \mu (dy)
\\
+ \lambda \int_E \{ a \varphi(x) - f(x) \} \varphi(x) \mu (dx)
\end{multline*}
for all $f$ and $a$.

\page{73}
Fix a measure $ \mu_0 $ on $E$ with $ \< \varphi, \mu_0 \> < \infty $.
Set 
$$
\mu^B_0 = 1_B \mu_0 , \quad \lambda^B_0 = \< \varphi 1_{B^c}, \mu_0
\>.
$$
Recall from \S2 that, for each compact set $B$, the equation
\begin{equation}
\eqlabel{LB} 
(\mu^B_t, \lambda^B_t) = (\mu^B_0, \lambda^B_0 ) + \int^t_0 L^B (\mu^B_s,
\lambda^B_s) \, ds
\end{equation}
has a unique solution, which is a continuous map
$$
t \mapsto (\mu^B_t, \lambda^B_t): [0, \infty) \rightarrow \cal M^+_B \times
\R^+.
$$
Consider \page{74} now a sequence of integer-valued  measures $ X^n_0 $. For
each $n$, denote by $ (X^n_t)_{t \geq 0} $ and $ (X^{B,n}_t,
\Lambda^{B,n}_t)_{t \geq 0} $ the stochastic coalescent and the
coupled family of approximations constructed above, starting from $
X^n_0 $. Set
$$
\eqalign{
\tilde X^n_t &= n^{-1} X^n_{n^{-1}t },\cr
(\tilde X^{B,n}_t, \tilde \Lambda^{B,n}_t) & = n ^{-1}
(X^{B,n}_{n^{-1}t}, \Lambda^{B,n}_{n^{-1}t}).
}
$$
We shall \page{75} need a mild continuity condition on $K$. Denote by $ S(K)
\subseteq E \times E $ the set of discontinuity points of $K$ and by $
\mu^{* n}_0 $ the $n$-th convolution power of $ \mu_0 $. Our
assumption is that
\begin{equation}
\eqlabel{CK}
(\mu^{*n}_0)^{\otimes 2} (S(K)) = 0, \quad \text{ for all } n \geq 1.
\end{equation}
This condition is verified, in particular, when $ S(K) $ has Lebesgue
measure zero and $ \mu_0 $ is absolutely continuous.

For the purposes of the next proposition, we also need an analogous
condition on the compact set $ B $:
\begin{equation}
\eqlabel{CB}
\mu^{* n}_0 (\partial B) = 0, \quad \text{ for all } n \geq 1.
\end{equation}
This condition is verified, for any given $ \mu_0 $, for all but
countably many closed intervals in $E$.
\page{76}
\begin{proposition}
\label{4.2}
Assume conditions {\rm \eqref{varphiK}}, {\rm \eqref{varphimu}}, 
{\rm \eqref{CK}}, {\rm \eqref{CB}}. Suppose that 
$$
d ( \tilde X^{B,n}_0, \mu^B_0) \rightarrow 0, \quad | \tilde
\Lambda^{B,n}_0 - \lambda^B_0 | \rightarrow 0
$$
as $ n \rightarrow \infty $. Then, for all $ t \geq 0 $
$$
\sup_{s \leq t} d(\tilde X^{B,n}_s, \mu^B_s) \rightarrow 0, \quad
\sup_{s \leq t} |\tilde \Lambda^{B,n}_s - \lambda^B_s | \rightarrow 0
$$
in probability as $ n \rightarrow \infty $.
\end{proposition}

\begin{proof}
Set $ \Lambda = \sup_n \< \varphi , \tilde X^n_0 \> $ and note that $
\Lambda < \infty $.
\page{77}
By rescaling \eqref{L} and \eqref{Q}, we see that, for all 
$B$, all bounded measurable functions $f$ and all $ a \in \R $
\begin{equation}
\eqlabel{Ln}
M^n_t = \sqrt{n}\left(\< f, \tilde X^{B,n}_t \> + a \tilde \Lambda^{B,n}_t
- \< f, \tilde X^{B,n}_0 \> - a \tilde \Lambda^{B,n}_0 - \int^t_0
L^{B,(n)} (\tilde X^{B,n}_s, \tilde \Lambda^{B,n}_s) (f,a) \, ds \right)
\end{equation}
is a martingale, having previsible increasing process
\begin{equation}
\eqlabel{Qn}
\< M^n \>_t = \int^t_0 Q^{B,(n)} (\tilde X^{B,n}_s, \tilde
\Lambda^{B,n}_s)(f,a) \, ds
\end{equation}
\page{78}
where
$$
\eqalign{
L^{B,(n)}(\mu, \lambda) & = n^{-2} L^{B,(1)} (n \mu, n \lambda), \cr
Q^{B,(n)} (\mu, \lambda) & = n^{-2} Q^{B,(1)} (n \mu, n \lambda).
}
$$
There is a constant $ C < \infty $, depending only on  $ B, \Lambda $
and $ \varphi $ such that
$$
\eqalign{ 
|L^B ( \tilde X^{B,n}_t, \tilde \Lambda^{B,n}_t)(f,a)| & \leq C(\|f\|
+ |a|),\cr
|(L^B - L^{B,(n)})(\tilde X^{B,n}_t, \tilde \Lambda^{B,n}_t)(f,a)| & \leq
C n^{-1} (\|f\| + |a|), \cr
|Q^{B,(n)} (\tilde X^{B,n}_t, \tilde \Lambda^{B,n}_t)(f,a)| & \leq
C(\|f\| + |a|)^2.
}
$$
Hence by the same argument as in Theorem \ref{4.1}, the laws of the
sequence $ (\tilde X^n, \tilde \Lambda^n) $ are tight in $ D([0,
\infty), \cal M_B \times \R) $. Indeed, similarly, the laws of the
sequence $ (\tilde X^n, \tilde \Lambda^n, I^n, J^n) $ 
are tight in $ D([0, \infty), \cal M_B \times \R \times \cal M_{B \times B}
\times \cal M_{B \times B}) $, where
$$
\eqalign{
I^n_t (dx, dy) &= K (x,y) 1_{x + y \in B} \tilde X^n_t (dx) \tilde
X^n_t(dy),\cr
J^n_t(dx, dy) &= K (x,y) 1_{x + y \not\in B} \tilde X^n_t (dx) \tilde
X^n_t (dy).
}
$$
\page{79}
Denote by $ (X,\Lambda, I, J) $ some weak limit point of this
sequence, which, by passing to a subsequence and the usual argument of
Skorohod, we may regard as a pointwise limit in $ D([0, \infty), \cal M_B
\times \R \times \cal M_{B \times B} \times \cal M_{B \times B}) $. Then there
exist bounded measurable functions
$$
I, J: \Omega \times [0, \infty) \times B \times B \rightarrow [0,
\infty)
$$
symmetric on $ B \times B $, such that
$$
\eqalign{
I_t(dx,dy) &= I(t,x,y) X_t (dx) X_t(dy),\cr
J_t (dx, dy) &= J(t,x,y) X_t(dx) X_t(dy)
}
$$
in $\cal M_{B\times B}$ and such that
$$
\eqalign{
I(t,x,y) & = K (x,y) 1_{x + y \in B} \cr
J(t,x,y) & = K (x,y) 1_{x + y \not\in B}.
}
$$
whenever $ (x,y) \not\in S(K) $ and $ x + y \not\in \partial B $.
Moreover we can pass to the limit in \eqref{Ln} to obtain, for all
continuous functions $f$ and all $ a \in \R $, for all $ t \geq 0 $,
almost surely
\page{80}
\begin{equation}
\eqlabel{E}
\hskip-1in
\eqalign{
\<(f,a), (X_t, \Lambda_t ) \> &= \<(f,a), (X_0, \Lambda_0) \> \cr
& \quad + \frac{1}{2} \int_0^t\int_{E \times E}\{ f(x + y) - f(x) - f(y) \} I(s,x,y) X_s(dx) X_s(dy)ds\cr
& \quad + \frac{1}{2} \int_0^t\int_{E \times E} \{ a \varphi (x + y) - f(x) - f(y) \}
J(s,x,y) X_s(dx) X_s(dy)ds \cr
& \quad + \int_0^t\Lambda_s \int_E \{ a \varphi(x)- f(x) \} \varphi(x) X_s (dx)ds
}
\hskip-1in
\end{equation}
By the remarks following the proof of Proposition \ref{2.2}, this
equation forces $ X_t \otimes X_t $ to be absolutely continuous with
respect to
$$
\sum^\infty_{n = 1} (\mu^{* n}_0)^{\otimes 2}.
$$
Hence by the assumptions \eqref{CK}, \eqref{CB}, we can replace  
$I(t,x,y)$ by $ K (x,y) 1_{x + y \in B} $  and $ J(t,x,y) $ by $K (x,y) 1_{x + y \not\in B}$ in \eqref{E}.
But this is now equation \eqref{LB} which has a unique solution $
(\mu^B_t , \lambda^B_t)_{t \geq 0}$.
\page{81}
We have shown that the unique weak limit point of $ (\tilde X^n,
\tilde \Lambda^n) $ in $ D ([0, \infty), \cal M_B) \times \R $ is the
continuous deterministic path $ (\mu^B_t, \lambda^B_t)_{t \geq 0} $,
which proves the proposition.
\end{proof}

\page{82}
We consider now the special case where $ \mu_0 $ is a probability measure
on $ \N = \{ 1,2, \ldots \} $. Here we can replace the method of weak
convergence in Proposition 4.2 by a
more direct approach using differential equations. The benefit in this
approach, besides greater transparency, is that we can establish a
rate of convergence, which is in principle computable. This would be
needed if one wished, in practice, to assess whether Smoluchowski's
equation provided a tolerable approximation to the stochastic
coalescent. Since the stochastic coalescent already makes a mean-field
approximation---it is assumed we neglect spatial variations in the
particle mass distribution---we are effectively assuming there is some
external spatial mixing. The relevant particle number $n$ is then the
number of particles in the largest region which is mixed to
equilibrium in unit time.
\page{83}
\begin{proposition}
\label{4.3}
Assume conditions {\rm \eqref{varphiK}} and {\rm \eqref{varphimu}}. Suppose that $
\mu_0 $ is supported on $ \N $. Let $B$ be a finite subset of $ \N $.
Then there is a constant $ C < \infty $, depending only on $K$,
$ \varphi, \mu_0 $, and $B$ such that, for all $ n \geq 1 $, for $
\tilde X^B_t = n^{-1} X^B_{n^{-1} t} $ and $ \tilde \Lambda^B_t =
n^{-1} \Lambda^B_{n^{-1}t} $, for all $ 0 \leq \delta \leq t $,
if $
\delta_0 = \|\tilde X^B_0 - \mu^B_0 \| + | \tilde \Lambda^B_0 -
\lambda^B_0|\leq 1 $ then
$$
\PP\left(\sup_{s \leq t} \{ \|\tilde X^B_s - \mu^B_s \| + | \tilde
\Lambda^B_s - \lambda^B_s | \} >  (\delta_0 + \delta) e^{Ct}\right) \leq
Ce^{-n \delta^2/Ct}.
$$
\page{84}
\end{proposition}
\begin{proof}
We regard $ (\tilde X^B_t, \tilde\Lambda^B_t)_{t \geq 0} $ as taking values
in the finite-dimensional vector space $ \R^B \times \R $. Recall from
the proof of Proposition \ref{4.2} that
\begin{equation}
\eqlabel{MD}
M_t = \sqrt{n} \{ (\tilde X^B_t, \tilde \Lambda^B_t) - (\tilde X^B_0,
\tilde \Lambda^B_0) - \int^t_0 L^{B,(n)} (\tilde X^B_s, \tilde
\Lambda^B_s)\,ds \}
\end{equation}
is a martingale, having previsible increasing process
$$
\< M \>_t = \int^t_0 Q^{B,(n)} (\tilde X^B_s, \tilde \Lambda^B_s) \,
ds.
$$
\page{85}

We recall a form of the exponential martingale inequality for
martingales $M$ whose jumps are bounded uniformly by $ A \in [0, \infty) $
and which have a continuous previsible increasing process $ \< M \>$: for all $ \theta \geq 0 $ we have
$$
\PP (\sup_t M_t \geq \delta \text{ and } \< M \>_\infty \leq \varepsilon)
\leq \exp \{ - \theta \delta + \frac{1}{2} \theta^2 e^{\theta A} \varepsilon \}.
$$
Let $ f : B \rightarrow \R $ and $ a \in \R $ be given, with $ \|f \|
+ |a| \leq 1 $. Consider the martingale
$$
M^{f,a}_t = \< (f,a),M_t\>.
$$
The jumps of $ M^{f,a} $ are bounded uniformly by $ Cn^{-\frac{1}{2}} $ for
some $ C < \infty $, depending on $ \varphi $ and $B$. The  process $
\< M^{f,a} \> $ is continuous and satisfies 
$$
\< M^{f,a}\>_t \leq Ct 
$$
for some $ C < \infty $, depending on $ \varphi, B $ and $ \mu_0 $.
So, by the exponential martingale inequality
$$
\PP \big(\sup_{s \leq t} M^{f,a}_s \geq \delta\big) \leq \exp \{ - \theta
\delta + \frac{1}{2} C \theta^2 t e^{C \,\theta/\sqrt{n}}\}.
$$
Assume that $ \delta \leq \sqrt{n}t $ and take $ \theta = \delta/(3Ct)
$. Then $ C \theta/\sqrt{n} \leq 1/3 $, so $ e^{C \theta/\sqrt{n}}
\leq 3/2 $ and so
$$
\PP \big(\sup_{s \leq t} M^{f,a}_s \geq \delta \big) \leq e^{-
\delta^2/4 Ct}.
$$
Since $B$ is finite, we deduce that for some $ C < \infty $, depending
on $ \varphi $, $B$ and $ \mu_0 $
$$
\PP \big(\sup_{s \leq t} \|M_s\| \geq \delta\big) \leq Ce^{-
\delta^2/Ct}
$$
whenever $ \delta \leq \sqrt{n} t $. Hence
$$
\PP \big(\sup_{s \leq t} \|n^{-\frac{1}{2}} M_s \| \geq \delta) \leq Ce^{-n
\delta^2/Ct} 
$$
whenever $ \delta \leq t $. Note also the estimate
$$
\|(L^B - L^{B,(n)})(\tilde X^B_t,\tilde \Lambda^B_t )\| \leq Cn^{-1},
\quad \text{ for all } t \geq 0,
$$
for some $ C < \infty $, depending on $ \varphi $ and $ B $.
\page{87}
Set $ Y_t = (\tilde X^B_t, \tilde \Lambda^B_t) - (\mu^B_t,
\lambda^B_t) $ and subtract the equations \eqref{MD}, \eqref{B} to
obtain
$$
Y_t = R_t + \int^t_0 L^B_s (Y_s) \, ds
$$
where
$$
R_t = Y_0 + n^{-\frac{1}{2}} M_t + \int^t_0 (L^B - L^{B,(n)})(\tilde X^B_s,
\tilde \Lambda^B_s) \, ds
$$
and where
$$
\eqalign{
\<(f,a), L^B_t (\mu, \lambda) \> 
& = \frac{1}{2} \int_{E \times E} \{ f(x + y) 1_{x + y \in B} + a \varphi(x +
y) 1_{ x + y \not\in B} - f(x) - f(y)\} \cr
& \qquad\qquad\times K (x,y) (\tilde X^B_t + \mu^B_t) (dx) \mu(dy)\cr
& \quad + (\tilde \Lambda^B_t + \lambda^B_t ) \int_E (a \varphi(x) - f(x))
\varphi(x) \mu(dx)\cr
&  \quad + \lambda \int_E (a \varphi(x) - f(x)) \varphi(x) (\tilde X^B_t +
\mu^B_t) (dx).
}
$$
Note the estimate
$$
\|L^B_t (\mu, \lambda)\| \leq C \|(\mu, \lambda)\|/2
$$
where $ C < \infty $ depends on $ \varphi $, $B$ and $ \mu_0 $. Set $
g(t) = \sup_{s \leq t}\|Y_s\| $ and $ r(t) = \sup_{s \leq t} \|R_s\|.
$ Then
$$
g(t) \leq r(t) + \frac{1}{2} C \int^t_0 g(s) \, ds
$$
\page{88}
so $ g(t) \leq r(t) e^{C t/2} $. Now, for $ \delta \leq t $, we have
$$
\PP(r(t) \geq g(0) + \delta/2 + Ct/n)\leq Ce^{-n \delta^2/Ct}.
$$
We may assume that $ C \geq 4 $. If $ \delta^2 \leq Ct/n $ we have
nothing to prove. Otherwise
$$
Ct/n < \delta^2 \leq \tfrac{1}{2} \delta e^{Ct/2}
$$
so
$$
\PP(g(t) \geq (\delta_0 + \delta) e^{Ct}) \leq \PP(r(t) \geq g(0) +
\delta/2 + Ct/n) \leq Ce^{-n \delta^2/Ct}.
\eqno\qed
$$
\let\qed\relax
\end{proof}
\page{89}
Here is the main result of this section.
\begin{theorem}
\label{4.4}
Let $ K: E \times E \rightarrow [0, \infty) $ be a symmetric
measurable function and let $ \mu_0 $ be a measure on $E$. Assume that
$$
(\mu^{* n}_0)^{\otimes 2} (S(K)) = 0 \quad \text{ for all } n \geq 1
$$
where $ S(K) $ denotes the discontinuity set of $K$. Assume also that,
for some continuous sublinear function $ \varphi: E \rightarrow (0,
\infty) $,
$$
K (x,y) \leq \varphi(x) \varphi(y) \quad \text{ for all } x, y \in E
$$
and that $ \< \varphi , \mu_0 \> < \infty $ and $ \< \varphi^2, \mu_0 \>
< \infty  $. Denote by $ (\mu_t)_{t < T} $ the maximal strong solution
provided by Theorem \ref{2.1}. Let $ (X^n_t)_{t \geq 0} $ be a
sequence of stochastic coalescents, with coagulation kernel $K$. Set $
\tilde X^n_t = n^{-1} X^n_{n^{-1}t} $ and suppose that
$$
d(\varphi \tilde X^n_0, \varphi \mu_0 ) \rightarrow 0
$$
as $ n \rightarrow \infty $. Then, for all $ t < T $
$$
\sup_{s \leq t} d (\varphi \tilde X^n_s, \varphi \mu_s) \rightarrow 0
$$
in probability, as $ n \rightarrow \infty $.
\page{90}

Moreover, if $ \mu_0 $ is supported on $ \N $, then for all $ t < T $
and all $ \delta >  0 $ there are constants $ \delta_0 >  0 $ and $  C
< \infty  $, depending only on $K$, $ \mu_0 $, $ \varphi $, $t$ and $
\delta $, such that, for all $n$,
$$
\|\varphi(\tilde X^n_0 - \mu_0)\| \leq \delta_0
$$
implies
$$
\PP\big(\sup_{s \leq t}\|\varphi(\tilde X^n_s - \mu_s)\| >
\delta\big) \leq e^{-n/C}.
$$
\page{91}
\end{theorem}
\begin{proof}
Fix $ \delta >  0 $ and $ t < T $. Since $ (\mu_t)_{t < T} $ is
strong, we can find a compact set $B$ satisfying \eqref{CB} and such
that $ \lambda^B_t < \delta/2 $. Now 
$$
d(\varphi \tilde X^n_0, \varphi \mu_0) \rightarrow 0
$$
so
$$
d(\tilde X^{B,n}_0, \mu^B_0) \rightarrow 0, \quad |\tilde
\Lambda^{B,n}_0 - \lambda^B_0 | \rightarrow 0.
$$
Hence, by Proposition \ref{4.3}
$$
\sup_{s \leq t} d (\tilde X^{B,n}_s, \mu^B_s) \rightarrow 0, \quad
\sup_{s \leq t} |\tilde \Lambda^{B,n}_s - \lambda^B_s| \rightarrow 0
$$
in probability as $ n \rightarrow \infty $. Since $ \{ \mu^B_s: s \leq
t \} $ is compact, we also have
$$
\sup_{s \leq t} d(\varphi  \tilde X^{B,n}_s, \varphi \mu^B_s)
\rightarrow 0
$$
in probability as $ n \rightarrow \infty $. By Proposition \ref{2.5}
and by \eqref{XI}, for $ s \leq t $
$$
\eqalign{
\|\varphi (\mu_s - \mu^B_s) \| & = \<  \varphi, \mu_s - \mu^B_s \>
\leq \lambda^B_s \leq \lambda^B_t < \delta/2 \cr
\|\varphi(\tilde X^n_s - \tilde X^{B,n}_s)\| & = \< \varphi, \tilde
X^n_s - \tilde X^{B,n}_s \> \leq \tilde \Lambda^{B,n}_s \leq \tilde
\Lambda^{B,n}_t \cr
& \leq \lambda^B_t + | \tilde \Lambda^{B,n}_t - \lambda^B_t | \cr
& \leq \delta/2 + |\tilde \Lambda^{B,n}_t - \lambda^B_t |.
}
$$
\page{92}
Now
$$
\eqalign{
d(\varphi \tilde X^n_s, \varphi \mu_s) &\leq \|\varphi(\tilde X^n_s -
\tilde X^{B,n}_s)\| + d(\varphi \tilde X^{B,n}_s, \varphi \mu^B_s) +
\|\varphi(\mu_s - \mu^B_s)\| \cr
& \leq \delta + d(\varphi \tilde X^{B,n}_s, \varphi \mu^B_s) + |\tilde
\Lambda^{B,n}_t - \lambda^B_t |
}
$$
so
$$
\PP (\sup_{s \leq t} d(\varphi \tilde X^n_s, \varphi \mu_s) >  \delta)
\rightarrow 0
$$
as $ n \rightarrow \infty $, as required.

In the case where $ \mu $ is supported on $ \N $, we can argue
similarly, replacing the weak metric $d$ by the total variation norm
and replacing Proposition \ref{4.3} by Proposition \ref{4.4}, to
arrive at the desired conclusion.
\end{proof}

\begin{corollary}
\label{4.5}
Let $K$, $ \mu_0 $, and $ \varphi $ be as in Theorem \ref{4.4}. Assume
in addition that $ \mu_0 $ is a probability measure and that $ \tilde
X^n_0 $ is the empirical distribution of a sample of size $n$ from $
\mu_0 $. Then, for all $ t < T $
$$
\sup_{s \leq t}\, d(\varphi \tilde X^n_s, \varphi \mu_s) \rightarrow 0
$$
in probability, as $ n \rightarrow \infty $. Moreover, if $ \mu_0 $ is
supported on $ \N $, and if $ \< e^{\alpha \varphi}, \mu_0 \> < \infty
$ for some $ \alpha >  0 $, then, for all $ t < T $ and all $ \delta >
0 $, there is a constant $ C < \infty $, depending only on $K$, $
\mu_0 $, $ \varphi $, $t$ and $ \delta $, such that for all $n$
$$
\PP(\sup_{s \leq t} \|\varphi(\tilde X^n_s - \mu_s) \| >  \delta) \leq
e^{-n/C}.
$$
\end{corollary}
\begin{proof}
For general $ \mu_0 $, it suffices to note that
$$
d (\varphi \tilde X^n_0 , \varphi \mu_0) \rightarrow 0
$$
almost surely as $ n \rightarrow \infty $, by the strong law of large
numbers, and to apply Theorem \ref{4.4}.
\page{94}

Suppose now that $ \mu_0 $ is supported on $ \N $. We have
$$
%\eqalign{
\|\varphi (\tilde X^n_0 
%&
- \mu_0) \|
%\cr
%&
\leq 2 \< \varphi 1_{(N, \infty)}, \mu_0 \> + | \< \varphi 1_{(N, \infty)}, \tilde X^n_0 - \mu_0 \> | 
%\cr
%&
+ \< \varphi 1_{(0,N]}, |\tilde X^n_0 - \mu_0 | \> 
%\cr & 
= I_1 + I_2 + I_3.
%}
$$
We can choose $N$ so that $ I_1 \leq \delta_0/2 $. Then by standard
exponential estimates we can find $ C < \infty $, depending on $
\mu_0, \varphi, N $ and $ \delta_0 $, such that
$$
\PP(I_2 + I_3 >  \delta_0/2) \leq e^{-n/C}.
$$
On combining this estimate with that found in Theorem \ref{4.4}, we
deduce
$$
\PP(\sup_{s \leq t} \|\varphi (\tilde X^n_s - \mu_s)\| >  \delta )
\leq 2e^{-n/C}
$$
as required.
\end{proof}
\marginpar{Author: add missing reference for \cite{Ch}.}

\bibliographystyle{alpha}
\bibliography{smol}

\begin{thebibliography}{McL64}

\bibitem[Ald]{A}
D.~J. Aldous.
\newblock Deterministic and stochastic models for coalescence (aggregation,
  coagulation):~a review of the mean-field theory for probabilists.
\newblock Preprint.
\newblock See www.stat.berkeley.edu/users/aldous.

\bibitem[BC90]{BC}
J.~M. Ball and J.~Carr.
\newblock The discrete coagulation-fragmentation equations: existence,
  uniqueness, and density conservation.
\newblock {\em J. Statist. Phys.}, 61(1-2):203--234, 1990.

\bibitem[Cha43]{Ch}
S.~Chandrasekhar.
\newblock Stochastic problems in physics and astronomy.
\newblock {\em Rev. Modern Physics}, 15:1--89, 1943.

\bibitem[CK]{CK}
J.~M.~C. Clark and V.~Katsouros.
\newblock Stable growth of a coarsening turbulent froth.
\newblock Preprint.

\bibitem[DS96]{DS}
P.~B. Dubovski{\u\i} and I.~W. Stewart.
\newblock Existence, uniqueness and mass conservation for the
  coagulation-fragmentation equation.
\newblock {\em Math. Methods Appl. Sci.}, 19(7):571--591, 1996.

\bibitem[EK86]{EK}
S.~N. Ethier and T.~K. Kurtz.
\newblock {\em Markov processes: characterization and convergence}.
\newblock Wiley Series in Probability and Mathematical Statistics. Wiley, New
  York, 1986.

\bibitem[Hei92]{H}
Ole~J. Heilmann.
\newblock Analytical solutions of {S}moluchowski's coagulation equation.
\newblock {\em J. Phys. A}, 25(13):3763--3771, 1992.

\bibitem[Jak86]{J}
Adam Jakubowski.
\newblock On the {S}korokhod topology.
\newblock {\em Ann. Inst. H. Poincar\'e Probab. Statist.}, 22(3):263--285,
  1986.

\bibitem[Jeo]{Je}
I.~Jeon.
\newblock Gelation phenomena.
\newblock Preprint.

\bibitem[Kur]{K}
T.~G. Kurtz.
\newblock Working paper on coalescence models.
\newblock Private communication.

\bibitem[Lus78]{Lu}
A.~A. Lushnikov.
\newblock Certain new aspects of the coagulation theory.
\newblock {\em Izv. Atmos. Ocean Phys.}, 14:738--743, 1978.

\bibitem[Mar68]{M}
Allan~H. Marcus.
\newblock Stochastic coalescence.
\newblock {\em Technometrics}, 10:133--143, 1968.

\bibitem[McL62]{McL1}
J.~B. McLeod.
\newblock On an infinite set of nonlinear differential equations.
\newblock {\em Quart. J. Math. Oxford}, 13:119--128, 1962.

\bibitem[McL64]{McL2}
J.~B. McLeod.
\newblock On the scalar transport equation.
\newblock {\em Proc. London Math. Soc.}, 14:445--458, 1964.

\bibitem[Pol84]{P}
David Pollard.
\newblock {\em Convergence of stochastic processes}.
\newblock Springer Series in Statistics. Springer-Verlag, New York, 1984.

\bibitem[vS16]{Sm}
M.~van Smoluchowski.
\newblock Drei {V}ortr\"age \"uber {D}iffusion, {B}rownsche {B}ewegung und
  {K}oagulation von {K}olloidteilchen.
\newblock {\em Physik. Z.}, 17:557--585, 1916.

\bibitem[Whi80]{W}
Warren~H. White.
\newblock A global existence theorem for {S}moluchowski's coagulation
  equations.
\newblock {\em Proc. Amer. Math. Soc.}, 80(2):273--276, 1980.

\end{thebibliography}
\end{document}